\documentclass[12pt,leqno]{amsart}  
\usepackage{amsmath,amstext,amsthm,amssymb,amsxtra}
\usepackage[colorlinks,citecolor=red,pagebackref,hypertexnames=false]{hyperref}
\usepackage[backrefs]{amsrefs}
\allowdisplaybreaks
\swapnumbers
\usepackage{pgf}
\numberwithin{equation}{section}

\allowdisplaybreaks
\swapnumbers
 
\theoremstyle{plain} 
\newtheorem{lemma}[equation]{Lemma} 
\newtheorem{proposition}[equation]{Proposition} 
\newtheorem{theorem}[equation]{Theorem} 
 
\newtheorem{conjecture}[equation]{Conjecture}
\newtheorem*{carleson}{Carleson's Theorem}

\theoremstyle{definition}

\theoremstyle{remark}
\newtheorem{remark}[equation]{Remark}
\newtheorem*{Acknowledgment}{Acknowledgment}

\def\t#1/{Theorem~\ref{t#1}}   \def\Tcap#1/{Theorem~\ref{t#1}} 
\def\c#1/{Corollary~\ref{c#1}}   \def\C#1/{Corollary~\ref{c#1}} 
\def\l#1/{Lemma~\ref{l#1}}        \def\Lcap#1/{Lemma~\ref{l#1}}  
\def\q#1/{Question~\ref{q#1}}	 \def\Pcap#1/{Proposition~\ref{p#1}}
\def\s#1/{Section~\ref{s#1}}      
\def\e#1/{(\ref{e#1})}
\def\r#1/{\ref{r#1}}
\def\f#1/{Figure~\ref{f#1}}

\def\Label #1 {\label{#1}}

\def\norm#1.#2.{\lVert#1\rVert_{#2}}
\def\Norm#1.#2.{\bigl\lVert#1\bigr\rVert_{#2}}
\def\NOrm#1.#2.{\Bigl\lVert#1\Bigr\rVert_{#2}}
\def\NORm#1.#2.{\biggl\lVert#1\biggr\rVert_{#2}}
\def\NORM#1.#2.{\Biggl\lVert#1\Biggr\rVert_{#2}}

\def\ip#1,#2.{\langle #1,#2\rangle}
\def\Ip#1,#2.{\bigl\langle#1,#2\bigr\rangle}
\def\IP#1,#2.{\Bigl\langle#1,#2\Bigr\rangle}

\def\abs#1{\lvert#1\rvert}
\def\Abs#1{\bigl\lvert#1\bigr\rvert}
\def\ABs#1{\Bigl\lvert#1\Bigr\rvert}

\def\XXint#1#2#3{{\setbox0=\hbox{$#1{#2#3}{\int}$}
     \vcenter{\hbox{$#2#3$}}\kern-.5\wd0}}

\def\ind#1 {{\mathbf 1}_{#1}}

\def\trans#1{\operatorname{Tr}_{#1}}
\def\modulate#1{\operatorname{Mod}_{#1}}
\def\dilate#1^#2{\operatorname{Dil}_{#1}^{#2}}

\def\ipf{\langle f,\varphi_s\rangle}

\def\dist #1,#2.{\operatorname{dist}(#1,#2)}
\def\dense #1{\operatorname{dense}(#1)}
\def\size #1{\operatorname{size}(#1)}
\def\COUnt #1{\operatorname{count}(#1)}
 \def\tree#1{$#1$tree} 
 \def\Tree{{\mathbf T}}
 
\def\mid{\,:\,}

\begin{document}

\title[Carleson's Theorem:Proof, Complements, Variations]{ Carleson's Theorem: \\
Proof, Complements, Variations }
 \author[M.T.~Lacey]{Michael T. Lacey
\\ { Georgia Institute of Technology}  }
\date{}

\thanks{This work has been
  supported by an NSF grant.} 

  \tableofcontents

\maketitle

 \section{Introduction}
  L.~Carleson's celebrated theorem of 1965 \cite{carleson} asserts the pointwise convergence of the
  partial Fourier sums of square integrable functions. We give  a  proof of this fact, in particular the proof of  
  Lacey and Thiele \cite{laceythielecarleson}, 
 as it can be presented in  brief self contained manner,  and a number of related results 
 can be seen by variants of the same argument.  We survey some of these variants, complements to Carleson's 
 theorem, as well as open problems.\footnote{This paper is an extended version of 
 the publication \cite {MR2091007}.}

 We are concerned with the Fourier transform on the real line, given by 
  \begin{equation*}
 \widehat f(\xi)=\int e^{-ix \xi} f(x)\;dx
  \end{equation*}
 for Schwartz functions $f$.  For such functions, it is an important elementary fact that one has Fourier inversion, 
  \begin{equation} \label{e.import-elem}
 f(x)=\lim_{N\to \infty}\frac1{2\pi}\int_{-N}^N \widehat f(\xi  )e^{ix \xi}\; d \xi,\qquad x\in\mathbb R,
 \end{equation} 
the inversion holding for   all Schwartz  functions $f$.   Indeed, 
 \begin{equation*}
\frac1{2\pi}\int_{-N}^N \widehat f(\xi  )e^{ix \xi}\; d \xi=D_N*f(x), 
 \end{equation*}
where $D_N(x):=\frac{\sin Nx }{\pi x} $ is the Dirchlet kernel.   

The convergence in (\ref{e.import-elem})  for Schwartz functions  follows from the classical facts 
 \begin{align*}
\int_{-\infty}^ \infty D_N(x)\; dx&{}=1,
\\
\lim_{N\to \infty} \int _{\abs x\ge \epsilon} D_N(x)\; dx&{}=0,\qquad \epsilon>0.
 \end{align*}

L.~Carleson's theorem asserts that (\ref{e.import-elem})   holds almost everywhere, for $f\in L^2(\mathbb R)$. 
The form of the Dirchlet kernel already points out the essential difficulties in establishing this theorem.  
That part of the kernel that is convolution with $\frac1 x $ corresponds to a singular integral. This 
can be done with the techniques associated to the Calder\'on Zygmund theory.  In addition, 
one must establish some uniform control for the oscillatory term $\sin N x $, which falls outside of 
what is commonly considered to be part of the Calder\'on Zygmund theory.  

For technical reasons, we find it 
easier to consider the equivalent one sided inversion, 
 \begin{equation}\label{e.Finversion} 
f(x)=\lim_{N\to \infty}\frac1{2\pi}\int_{-\infty}^N \widehat f(\xi  )e^{ix \xi}\; d \xi.
 \end{equation}

 Schwartz functions being dense in $L^2$, 
one need only show that  the set of functions for which a.e.~convergence holds is closed.  The standard method for doing so 
is to consider the maximal function below, which we refer to as the Carleson operator 
 \begin{equation}\label{e.carlesonoperator} 
\mathcal Cf(x):=\sup_N \ABs{ \int_{-\infty}^N \widehat f(\xi  )e^{ix \xi}\; d \xi},\qquad x\in\mathbb R.
 \end{equation}
There is a straight forward proposition.

 \begin{proposition}\label{p.maxclosed}  
Suppose that the Carleson operator satisfies 
 \begin{equation}\label{e.carweakL2}
\abs{ \{\mathcal C f(x)>\lambda\}}\lesssim{} \lambda^{-2}\lVert f\rVert_2^2,\qquad f\in L^2(\mathbb R),\ \lambda>0.
 \end{equation}
Then, the set of functions $f\in L^2(\mathbb R)$ for which (\ref{e.Finversion})  holds is closed and hence all of $L^2(\mathbb R)$. 
 \end{proposition}  

\begin{proof}  For $f\in{}L^2(\mathbb R)$, we should see that 
 \begin{equation*}
L_f:=\limsup_{N\to\infty}\Abs{ f(x)-\frac1{2\pi}\int_{-\infty}^N \widehat f (\xi)e^{ix\xi}\; d\xi }=0\quad\text{a.e.}
 \end{equation*}
To do so, we show that for all $\epsilon>0$,  $\abs{\{L_f>\epsilon\}}\lesssim\epsilon$.  We take $g$ to be a smooth compactly supported 
function so that $\norm f-g.2.\le\epsilon^{3/2}$.  Now 
Fourier inversion holds for $g$, whence  
$
L_f\le{}\mathcal C(f-g)+\abs{f-g}
$.  Then, by the weak type inequality, (\ref{e.carweakL2}) , we have 
 \begin{equation*}
\abs{ \{ \mathcal C(f-g)>\epsilon\}}\lesssim\epsilon^{-2}\lVert f-g\rVert_2^2\lesssim\epsilon.
 \end{equation*}
\end{proof}

This is a standard proposition, which holds in a general context, and serves as one of the prime motivations for 
considering maximal operators.  Note in particular that we are not at this moment claiming that $\mathcal C$ is a bounded 
operator on $L^2$.  Inequality (\ref{e.carweakL2})  is the so called weak $L^2$ bound, and we shall utilize the form of this
bound in a very particular way in the proof below.

 It was one of L.~Carleson's great achievements to invent a method to prove this weak type estimate.  
 
  \begin{carleson}\label{t.carleson} The estimate (\ref{e.carweakL2})  holds.  As a consequence, (\ref{e.Finversion})  holds for all $f\in L^2(\mathbb R)$, 
 for almost every $x\in\mathbb R$. 
  \end{carleson} 
 
 Carleson's original proof \cite{carleson} was extended to $L^p$, $1<p<\infty$, by R.~Hunt.  Also see \cite{MR39:3222}. 
 Charles Fefferman
 \cite{fefferman} 
 gave an alternate proof
that was influential by the explicit nature of it's   ``time frequency'' analysis, of which we have more more to say below.  
We  follow the proof of M.~Lacey and C.~Thiele \cite{laceythielecarleson}. More detailed comments on the
history of the proof, and related results will come later. 

The proof will have three  stages, the first being an appropriate decomposition of the Carleson operator. 
The second being an introduction of three Lemmas, which can be efficiently combined to give the proof of our Theorem.
The third being a proof of the 
Lemmas. 

\bigskip 
We do not keep track of the value of generic absolute constants, instead using the  notation $A\lesssim{}B$ iff $A\le{}KB$ 
for some constant $K$.   And $A\simeq B$ iff $A\lesssim{}B$ and $B\lesssim{}A$.  The notation
$\ind A $ denotes the indicator function of the set $A$.
For an operator $T$, $\norm T.p.$ denotes the norm of $T$ as an operator from $L^p $ to itself.  

\begin{Acknowledgment}
These notes are based on a series of lectures given at the Erwin Schr\"odinger  
  Institute, in Vienna Austria. I an                                                                               
  indebted to the Institute for the opportunity to present these lectures.
  \end{Acknowledgment}

\section{Decomposition}

The Fourier transform is a constant times a unitary operator on $L^2(\mathbb R)$.  In particular, we shall take the Plancherel's identity for granted. 

 \begin{proposition}\label{p.planch}  For all $f,g\in L^2(\mathbb R)$,  
 \begin{equation*}
\ip f,g.=c\ip \widehat f,\widehat g.
 \end{equation*}
for appropriate constant $c=\frac1{2\pi}$.
 \end{proposition} 

The convolution of $f$ and $\psi$ is given by $f*\psi(x)=\int f(x-y)\psi(y)\; dy$.  We shall 
also assume the following Lemma.

 \begin{lemma}\label{l.convolve} If a bounded  linear operator $T$ on $L^2(\mathbb R)$ commutes with translations, then 
$Tf=\psi*f$, where $\psi$ is a distribution, which is to say a linear functional on Schwartz functions. 
In addition, the Fourier transform of $Tf$ is given by 
 \begin{equation*}
\widehat{Tf}=\widehat{\psi}\widehat f .
 \end{equation*}
 \end{lemma}

Let us introduce the operators associated to translation, modulation and dilation on the real line. 
 \begin{align}\label{e.translation} 
\trans y f(x)&{}:=f(x-y),
\\  \label{e.modulate}
\modulate \xi f(x)&{}:=e^{i \xi x}f(x),
\\ \label{e.dilate}
\dilate \lambda^p f(x)&{}:= \lambda^{-1/p}f(x/\lambda),\qquad 0<{}p\le\infty,\ \lambda>0.
 \end{align}
Note that the dilation operator preserves $L^p$ norm. 
These operators are  related through the Fourier transform, by 
 \begin{equation}\label{e.fouriertransformoperator}
 \widehat{\trans y}=\modulate {-y},\qquad \widehat{\modulate \xi}=\trans \xi,\qquad \widehat{\dilate \lambda^2}=\dilate {1/\lambda}^2
  \end{equation}
And we should also observe that the Carleson operator commutes with translation and dilation operators, while being invariant under modulation operators.  For any $y, \xi\in\mathbb R$, and $\lambda>0$,  
 \begin{equation*}
\trans y \circ\, \mathcal C=\mathcal C\circ \trans,\qquad \dilate \lambda^2 \circ\, \mathcal C=\mathcal C\circ \dilate \lambda^2,\qquad \mathcal C\circ\modulate \xi=\mathcal C.
  \end{equation*}
Thus, our mode of analysis should exhibit the same invariance properties.

We have phrased the Carleson operator in terms of  modulations of 
the operator $\operatorname P_-f(x)=\int_{-\infty}^0 \widehat f(\xi)e^{ix \xi}\; d \xi$, 
which is the Fourier projection on to negative frequencies.  
Specifically, since multiplication of $f$ by an exponential is associated with a translation of $\widehat f$, we have 
 \begin{equation}\label{e.carelsonP-}
\mathcal Cf=\sup_N \abs{\operatorname P_-(e^{iN\cdot}f)}
 \end{equation}

\bigskip 
A characterization of the operator $\operatorname P_-$ will be useful to us.

 \begin{proposition}\label{p.P-}  Up to a constant multiple, $\operatorname P_-$ is the unique bounded 
operator on $L^2(\mathbb R)$ which  (a) commutes with translation (b) commutes with dilations (c) has as it's kernel precisely those functions with frequency support on the positive axis.
 \end{proposition} 

\begin{proof}  Let $T$ be a bounded operator on $L^2(\mathbb R)$  which satisfies these three properties. 
Condition (a) implies that $T$ is given by  convolution with respect to a distribution.  Such operators are equivalently characterized in frequency variables by 
$\widehat{Tf}=\tau\widehat f$ for some bounded function $\tau$.   Condition (b) then implies that $\tau( \xi)=\tau(\xi/\abs\xi)$
for all $\xi\not=0$.  
A function $f$ is in the kernel of $T$ iff $\widehat f$ is supported on the zero set of $\tau$.  
Thus (c) implies that $\tau$ is identically $0$ on the positive real axis, and 
non--zero on the negative axis.  Thus, $T$ must be a multiple of $\operatorname P_-$.
\end{proof}

We move towards the tool that will permit us to decompose  Carleson operator, and take advantage of some combinatorics of the 
time--frequency plane.  
We let $\mathcal D$ be a choice of dyadic grids on the line. Of the different choices we can make, we take the grid to be one that 
is preserved under dilations by powers of $2$. That is 
 \begin{equation}\label{e.Ddef}
\mathcal D=\{ [j2^k,(j+1)2^k)\mid j,k\in\mathbb Z\}
 \end{equation}
Thus, for each interval $I\in\mathcal D$ and $k\in\mathbb Z$,  the interval $2^kI=\{2^kx\mid x\in I\}$ is also in $\mathcal D$.

  A {\em tile} is a rectangle $s\in\mathcal D\times \mathcal D$ that has area one. 
We write a tile as $s=I\times \omega$, thinking of the first interval as a time interval and the second as frequency.  The requirement 
of having area one is suggested by the uncertainty principle of the Fourier transform, or alternatively, our calculation of the 
Fourier transform of the dilation operators in (\ref{e.fouriertransformoperator}) . Let $\mathcal T$ denote the set of all tiles.  
While tiles all have area one, the ratio between the time and frequency coordinates is permitted to be 
arbitrary.  See \f.grids/ for a few possible choices of this ratio.

\begin{figure}
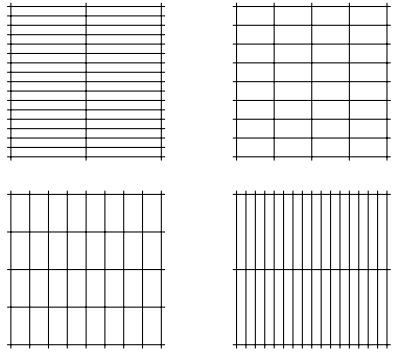
 

\begin{center}
\begin{pgfpicture} {0cm}{0cm}{4cm}{4.5cm}

\pgfgrid[step={\pgfpoint{1cm}{.125cm}}]{\pgfxy(-.05,2.45)}{\pgfxy(2.05,4.55)}

\pgfgrid[step={\pgfpoint{.5cm}{.25cm}}]{\pgfxy(2.95,2.45)}{\pgfxy(5.05,4.55)} 

\pgfgrid[step={\pgfpoint{.25cm}{.5cm}}]{\pgfxy(-.05,-.05)}{\pgfxy(2.05,2.05)}

\pgfgrid[step={\pgfpoint{.125cm}{1cm}}]{\pgfxy(2.95,-.05)}{\pgfxy(5.05,2.05)}

\end{pgfpicture}
\end{center}
\label{f.grids} 
 \caption{Four different aspects ratios for  tiles. Each fixed ratio gives rise to a tiling for the time frequency plane.} 
\end{figure}


Each dyadic interval is a union of its left and right halves, which are also dyadic.  For an interval $\omega$ we denote these as 
$\omega_-$ and $\omega_+$ respectively.  We are in the habit of associating frequency intervals with the vertical axis. 
So  $\omega_-$ will lie below   $\omega_+$. Associate to a  tile $s=I_s\times \omega_s$ the rectangles $s_{\pm}=I_s\times 
\omega_{s\pm}$. 
These two rectangles play complementary roles in our proof.  

Fix a Schwartz function $\varphi$ with $\ind {[-1/9,1/9]} \le\widehat \varphi\le\ind {[-1/8,1/8]} $.
Define a function associated to a tile $s$ by 
 \begin{equation}\label{e.zvf}
\varphi_s=\modulate {c(\omega_{s-})} \trans {c(I_s)} \dilate {\abs{I_s}}^2  \varphi
 \end{equation}
where $c(J)$ is the center of the interval $J$.  Notice that $\varphi_s$  has Fourier transform supported on $\omega_{s-}$,
and is highly localized in time variables around the interval $I_s$.  That is, $\varphi_s$ 
is essentially supported in the time--frequency plane on the rectangle $I_s\times \omega_{s-}$. 
Notice that the set of functions $\{\varphi_s\mid s\in\mathcal T\}$  has a set of invariances with respect to 
translation, modulation, and dilation that mimics those of the Carleson operator. 

\medskip 

It is our purpose to devise a decomposition of the projection $\operatorname P_-$ in terms of the tiles just introduced.  
To this end, for a choice of $\xi\in\mathbb R$, let 
 \begin{equation}\label{e.Qdef} 
\operatorname Q_ \xi f=\sum_{\substack{s\in\mathcal T}} \ind {\omega_{s+}}  (\xi)\ipf \varphi_s .
 \end{equation}
We should consider general values of $\xi$ for the reason that the dyadic grid distinguishes certain points as  being 
interior, or a boundary point, to an infinite chain of dyadic intervals.   And moreover, for a given $\xi$, only certain 
tiles can contribute to the sum above, those tiles being determined by the expansion of $\xi $ in a binary digits.  See 
\f.Q/.  Let us list  some  relevant properties of these operators.

\begin{figure}
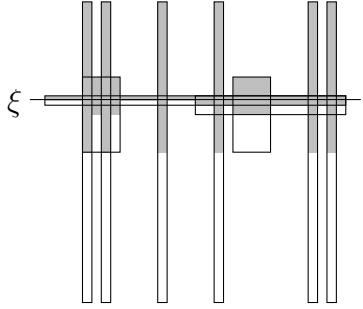
 

\begin{center} 
\begin{pgfpicture}{0cm}{0cm}{4cm}{4.5cm} 

\color{lightgray}
\begin{pgftranslate}{\pgfpoint{0cm}{0.5cm}}

\pgfrect[fill]{\pgfxy(0,1.9375)}{\pgfxy(4,.0625)}
\pgfrect[fill]{\pgfxy(2,1.875)}{\pgfxy(2,.125)}
\pgfrect[fill]{\pgfxy(.5,1.75)}{\pgfxy(.5,.5)}
\pgfrect[fill]{\pgfxy(2.5,1.75)}{\pgfxy(.5,.5)}

\pgfrect[fill]{\pgfxy(.5,1.25)}{\pgfxy(.125,2)}
\pgfrect[fill]{\pgfxy(.75,1.25)}{\pgfxy(.125,2)}
\pgfrect[fill]{\pgfxy(1.5,1.25)}{\pgfxy(.125,2)}
\pgfrect[fill]{\pgfxy(2.25,1.25)}{\pgfxy(.125,2)}
\pgfrect[fill]{\pgfxy(3.5,1.25)}{\pgfxy(.125,2)}
\pgfrect[fill]{\pgfxy(3.75,1.25)}{\pgfxy(.125,2)}

\color{black} 

\pgfputat{\pgfxy(-.4,1.9)}{\pgfbox[center,center]{$\xi$}} 
\pgfline{\pgfxy(-.2,1.95)}{\pgfxy(4.2,1.95)}  
\pgfrect[stroke]{\pgfxy(0,1.875)}{\pgfxy(4,.125)}
\pgfrect[stroke]{\pgfxy(2,1.75)}{\pgfxy(2,.25)}
\pgfrect[stroke]{\pgfxy(.5,1.25)}{\pgfxy(.5,1)}
\pgfrect[stroke]{\pgfxy(2.5,1.25)}{\pgfxy(.5,1)}
\pgfrect[stroke]{\pgfxy(.5,-.75)}{\pgfxy(.125,4)}
\pgfrect[stroke]{\pgfxy(.75,-.75)}{\pgfxy(.125,4)}
\pgfrect[stroke]{\pgfxy(1.5,-.75)}{\pgfxy(.125,4)}
\pgfrect[stroke]{\pgfxy(2.25,-.75)}{\pgfxy(.125,4)}
\pgfrect[stroke]{\pgfxy(3.5,-.75)}{\pgfxy(.125,4)}
\pgfrect[stroke]{\pgfxy(3.75,-.75)}{\pgfxy(.125,4)}

\end{pgftranslate}

\end{pgfpicture} \end{center} 
		\label{f.Q}  \caption{Some of the tiles that contribute to the sum for $\operatorname Q_\xi $.  
		The shaded areas are the tiles $I_s\times \omega_{s+} $.}
\end{figure}

 \begin{proposition}\label{p.Qprop}  For any $\xi$, the operator $\operatorname Q_ \xi$ is a bounded operator on $L^2$, with bound independent of $ \xi$.  Its kernel 
contains those functions with Fourier transform supported on $[\xi,\infty)$, and it is positive semidefinite.  Moreover, 
for each integer $k$, 
 \begin{gather}\label{e.Qdilate}
\operatorname Q_ \xi=\dilate {2^{-k}}^2 \operatorname Q_ {\xi2^{-k}} \dilate {2^{k}}^2
 \\ \label{e.Qtrans}
\operatorname Q_ {\xi,k}\trans {2^k}=\trans {2^k}Q _ {\xi,k},
 \end{gather}
where $\operatorname Q_{ \xi,k}=\sum_{\substack{s\in\mathcal T\\ \abs{I_s}\le{}2^k}} \ind {\omega_{s+}}  (\xi)\ipf \varphi_s $.
 \end{proposition} 

\begin{proof}  Let $\{\omega(n)\mid n\in\mathbb Z\}$ be the set of dyadic intervals for which  $\xi\in\omega(n)_+$, listed in increasing order, thus 
$\cdots \subset\omega(n)\subset \omega(n+1)\subset\cdots$.    Let $\mathcal T(n)=\{s\in\mathcal T\mid \omega_s=\omega(n)\}$, 
and 
 \begin{equation*}
\operatorname Q_{(n)}f=\sum_{s\in\mathcal T(n)} \ipf \varphi_s .
 \end{equation*}
The intervals $\omega(n)_-$ are disjoint in $n$, and since $\varphi_s$ has frequency support in $\omega_{s-}$, it follows that the operators $\operatorname Q_{(n)}$ are orthogonal in $n$.  The boundedness of $\operatorname Q_ \xi$ reduces therefore to the
 uniform boundedness of $\operatorname Q_{
(n)}$ in $n$.  

Two operators $\operatorname Q_{(n)}$ and $\operatorname Q_{(n')}$ differ by composition with a modulation operator and a dilation operator that preserves $L^2$ norms. 
Thus, it suffices to consider the $L^2$ norm bound of a  $\operatorname Q_{(n)}$ with  $\abs{I_s}=1$ for all $s\in\mathcal T(n)$.  Using the fact that $\varphi_s$ is a rapidly decreasing function, we see that 
 \begin{equation} \label{e.compare}
\abs{\ip \varphi_s, \varphi_{s'}.}\lesssim{} {\dist I_s,I_{s'}.}^{-4}.
 \end{equation}
Now that  the spatial length of the tiles is one, the 
 tiles are separated by integral distances.  Since $\sum_n n^{-4}<\infty$, 
 \begin{align*}
\lVert \operatorname Q_{(n)}f\rVert_2^2={}&\sum_{s\in\mathcal T(n)}\sum_{s'\in\mathcal T(n)}\ipf \ip \varphi_s, \varphi_{s'}. 
		\ip \varphi_{s'},f.
	\\
	{}\lesssim{}&\sup_{n\in\mathbb Z} \sum_{s\in\mathcal T(n)}\abs{ \ipf \ip f,\varphi_{(I_{s}+n)\times \omega(n)}.}
	\\
	{}\lesssim{}&\sum_{s\in\mathcal T(n)}\abs{ \ipf }^2.
 \end{align*}
The last inequality following by Cauchy-Schwarz.  The last sum is easily controlled, by simply bringing in the absolute values. 
Since $\abs{\ipf}^2\lesssim{}\int \abs{f}^2\abs{\varphi_s}\; dx$
 \begin{align*}
\sum_{s\in\mathcal T(n)}\abs{ \ipf }^2\le{}& \int \abs {f}^2\sum_{s\in\mathcal T(n)}\abs{\varphi_s}\; dx
  \\
  {}\le{}&\lVert f \rVert_2^2\sup_x \sum_{s\in\mathcal T(n)}\abs{\varphi_s(x)}
  \\
  {}\lesssim{}&\lVert f \rVert_2^2
  \end{align*}
This completes the proof of the uniform boundedness of $\operatorname Q_ \xi$.  

\medskip

Since all of the functions $\varphi_s$ that contribute to the definition of $\operatorname Q_ \xi$ have frequency support below $\xi$, the conclusion abut the kernel of the operator is obvious. And that it is positive semidefinite, observe that 
 \begin{equation}
\ip \operatorname Q_ \xi f,f.=\sum_{\substack{s\in\mathcal T\\ \xi\in \omega_{s+}}} \abs{\ipf }^2\ge0.
 \end{equation}
In particular, $\ip \operatorname Q_ \xi \varphi_s,\varphi_s.\not=0$ for $s\in\mathcal T(n)$. 

\medskip 
To see (\ref{e.Qdilate})  recall (\ref{e.fouriertransformoperator})  and our specific choice of grids, (\ref{e.Ddef}) .  To see (\ref{e.Qtrans}) , observe that if $I\in\mathcal D$ has length at most $2^k$, then $I+2^k$ is also in $\mathcal D$.

\end{proof}

As the lemma makes clear, $ \modulate {-\xi} \operatorname Q_ \xi\modulate \xi$ serves as an approximation to $\operatorname P_-$.    A limiting procedure will recover $\operatorname P_-$ exactly. Consider 
 \begin{equation} \label{e.Qaverage}
Q=\lim_{Y\to \infty} \int_{B(Y)} \dilate {2^{-\lambda}}^2 \trans {-y} \modulate {-\xi} \operatorname Q_ \xi 
\modulate {\xi}\trans {y}  \dilate {2^{\lambda}}^2 \mu(d \lambda,dy,d \xi).
 \end{equation}
Here, $B(Y)$ is the set $[1,2]\times[0,Y]\times [0,Y]$, and $\mu$ is normalized Lebesgue measure.  Notice that the dilations are 
given in terms of $2^\lambda$, so that in that parameter, we are performing an average with respect
to the multiplicative Haar measure on $\mathbb R_+$.

Apply the right hand side to a Schwartz function $f$.  It is easy to see that as $k\to-\infty$, the terms 
 \begin{equation*}
\modulate {\xi}\trans {y}  \dilate {2^{\lambda}}^2\operatorname Q_{\xi, k}f
 \end{equation*}
 tend to  zero uniformly in the parameters $\xi$, $y$, and $\lambda$, with a rate that depends upon $f$.
  Here, $\operatorname Q_{\xi,k}$ is as in (\ref{e.Qtrans}) . Similarly, 
as $k\to\infty$, the terms 
   \begin{equation*}
\modulate {\xi}\trans {y}  \dilate {2^{\lambda}}^2(\operatorname Q-\operatorname Q_{\xi, k})f
 \end{equation*}
also tend to zero uniformly.  Hence, 
the limit is seen to exist for all Schwartz functions.  By \Pcap.Qprop/, it follows that $Q$ is a bounded operator on $L^2$. 
That $Q$ is translation and dilation invariant follows from (\ref{e.Qtrans})  and (\ref{e.Qdilate}) . 
Its kernel contains those functions with Fourier transform supported on $[0,\infty)$. 
Finally, if we verify that $Q$ is not identically zero, we can conclude that it is a multiple of $\operatorname P_-$.  
But, for e.g.~$f=\modulate {-1/8}\varphi$, it is easy to see that 
 \begin{equation*}
\ip \operatorname Q_ \xi 
\modulate {\xi}\trans {y}  \dilate {2^{\lambda}}^2 f,\modulate {\xi}\trans {y}  \dilate {2^{\lambda}}^2 f .>0
 \end{equation*}
so that $Qf\not=0$.  Thus, $Q$ is a multiple of $\operatorname P_-$.

\medskip 

We can return to the Carleson operator.  An important viewpoint emphasized by C.~Fefferman's proof 
\cite{fefferman} is that we should linearize the supremum.  That is we consider a measurable map $N\mid \mathbb R\mapsto\mathbb R$, 
which specifies the value of $N$ at which the supremum in 
(\ref{e.carlesonoperator})  occurs.  Then, it suffices to bound the operator norm of the linear (not sublinear) operator
 \begin{equation*}
\operatorname P_-\modulate {N(x)}
 \end{equation*}
Considering (\ref{e.Qaverage}) , we set 
 \begin{equation*}
\mathcal C_N f(x)=\sum_{s\in\mathcal T}\ind \omega_{s+}   (N(x))   \ipf \varphi_s (x).
 \end{equation*}
Our main Lemma is then 

 \begin{lemma}\label{l.model}   There is an absolute constant $K$ so that for all measurable functions $N\mid \mathbb R\mapsto\mathbb R$, we 
have the weak type inequality 
 \begin{equation} \label{e.modelweaktype} 
\abs{ \{\mathcal C_N f>\lambda\}}\lesssim{} \lambda^{-2}\lVert f\rVert_2^2,\qquad \lambda>0,\ f\in L^2(\mathbb R).
 \end{equation}
 \end{lemma}

By the convexity of the weak $L^2$ norm, this theorem immediately implies the same estimate for $\operatorname P_-\modulate {N(x)}$, and so proves 
\t.carleson/. 
The proof of the lemma is obtained by combining the three estimates detailed in the next section.

\subsection{Complements}  
At the conclusion to the different sections of the proof, some complements to the ideas and techniques of the 
previous sections will be mentioned, but not proved.  These items can be considered as exercises. 

 \begin{remark}\label{r.avgconvolve}  For Schwartz functions $\varphi$ and $\psi$, set 
 \begin{gather*}
\operatorname A f:=\sum_{n\in\mathbb Z} \ip f,\trans n \varphi . \trans n \psi 
\\
\operatorname B f:=\int_0^1 \trans {-y} \operatorname A \trans y \; dy
 \end{gather*}
Then, $\operatorname B$  is a convolution operator, that is $\operatorname  Bf=\Psi*f$ for some function $\Psi$, which can be explicitly computed.    
 \end{remark}

 \begin{remark}\label{r.identity}  The identity operator is, up to a constant multiple, the 
unique bounded operator $\operatorname A$ on $L^2$ which commutes with all translation and modulation 
operators.  That is  $\operatorname A \mid L^2 \mapsto L^2 $, and for all $y\in\mathbb R$ and $\xi\in\mathbb R$, 
 \begin{equation*}
\trans{y}\operatorname A=\operatorname A\trans y,\qquad \modulate \xi \operatorname A=\operatorname A\modulate \xi .
 \end{equation*}
 \end{remark}

 \begin{remark}\label{r.Aj}   The operators  
 \begin{equation*}
\operatorname A_jf:=\sum_{\substack{ s\in \mathcal T\\ \abs{I_s}=2^j}} \ip f,\varphi_s . \varphi_s 
 \end{equation*}
are uniformly bounded operators on $L^2(\mathbb R)$, assuming that $\varphi $ is a Schwartz function.  
 \end{remark}

 \begin{remark}\label{r.avgident}  Assuming that $\varphi\not\equiv0 $, the operator below is a 
non--zero multiple of the identity operator on $L^2$. 
 \begin{equation*}
\int_0^{2^j} \int_0^{2^{-j}} \trans {-y} \modulate {-\xi} \operatorname  A_j \modulate \xi \trans y \; d\xi dy.
 \end{equation*}
 \end{remark}

\section{The Central Lemmas} 

Observe that the weak type estimate of \Lcap.model/  is implied by 
 \begin{equation} \label{e.WKT}
\abs{ \ip \mathcal C_N f,\ind E  .}\lesssim{}\norm f.2. \abs{E}^{1/2},
 \end{equation}
for all functions $f$ and sets $E$ of finite measure.  (In fact this inequality is 
equivalent to the weak $L^2$ bound.)

By linearity of $\mathcal C_N$, we may assume that $\norm f.2.=1$.  By the invariance of $\mathcal C_N$ under dilations by a factor of $2$ (with a change of measurable $N(x)$), we can assume that $1/2<\abs E\le1$.   Set 
 \begin{equation}\label{e.zf}
\phi_s=(\ind {\omega_{s+}} \circ  N)\varphi_s.
 \end{equation}
We shall show that 
 \begin{equation}\label{e.weaktypedual} 
\sum_{s\in\mathcal T}\abs{\ipf \ip \phi_s,\ind E .}\lesssim{}1.
 \end{equation}
To help keep the notation straight, note that $\varphi_s $ is a smooth function, adapted to the tile. On the other hand, 
$\phi_s $ is the rough function paired with the indicator set $\ind {\omega_{s+}} \circ N $.  
From this point forward, the function $f$ and the set $E$  are fixed.  We use data about these two objects to organize our proof. 

\smallskip

 As the sum above is over strictly positive quantities,
 we may consider all sums to be taken over some finite subset of tiles.  Thus, there is never any
 question that the sums we treat are 
 finite, and 
the iterative procedures we describe will all terminate.  The estimates we obtain will be independent of the 
fact that the sum is formally over a set of finite tiles. 

\smallskip 

We need some concepts to phrase the proof.  There is a natural partial order on tiles. Say that $s<s'$ iff $\omega_{s}\supset 
\omega_{s'}$ and $I_s\subset I_{s'}$.  Note that the time variable of $s$ is localized to that of $s'$, and the frequency 
variable of $s$ is 
similarly localized, up to the variability allowed by the uncertainty principle. Note that two tiles are incomparable 
with respect to the 
`$<$' partial order iff the tiles, as rectangles in the time frequency plane, do not intersect.   A ``maximal tile''{} will 
one that is maximal 
with respect to this partial order.  See figure 1.

We call a set of tiles $\Tree \subset\mathcal S$ a {\em tree} if there is a tile $I_\Tree\times \omega_{\Tree}$, called the
{\em top of the tree,} 
such that for all $s\in\Tree$, $s<I_\Tree\times \omega_{\Tree}$.  We note that the top is not uniquely defined.  An important point is that a tree top specifies a location in time variable for the tiles in the tree, namely inside $I_\Tree$, and localizes 
the  frequency variables,  identifying $\omega_{\Tree}$ as 
a nominal origin.

We say that {\em $\mathcal S$ has count at most $A$}, and write 
 \begin{equation*}
\COUnt {\mathcal S}<A
 \end{equation*}
iff $\mathcal S$ is a union   $\bigcup_{\Tree\in\mathcal T}\Tree$, where each $\Tree\in\mathcal T$ is a tree, and  \begin{equation*}
\sum_{\Tree\in\mathcal T}\abs{I_\Tree}<A.
 \end{equation*}

Fix $\chi(x)=(1+\abs x)^{-\kappa}$, where $\kappa$ is a large constant, whose exact value is unimportant to us.  Define  
 \begin{gather} \label{e.chiI} 
\chi_I:=\trans {c(I)}\dilate {\abs I}^1  \chi,
\\
\label{e.dense} 
\dense s:=\sup_{s<s'}\int_{ N^{-1}(\omega_{s'})} \chi_{I_{s'}}\; dx,
\\\nonumber 
\dense {\mathcal S}:=\sup_{s\in\mathcal S}\dense s,\qquad \mathcal S\subset\mathcal T.
 \end{gather}
The first and most natural definition of a ``density'' of a tile, would be ${\abs {I_s}}^{-1}\abs{N^{-1}(\omega_{s+})\cap I_s}$.  But 
$\varphi$ is supported on the whole real line, though does decay faster than any inverse of a  polynomial. 
We refer to this as a ``Schwartz tails problem.''
  The definition of density as  $\int_{N^{-1}(\omega_{s})} \chi_{I_s}\;dx $, as it turns out, is still not adequate. 
 That we should take the supremum over  $s<s'$ only becomes evident in the proof of the ``Tree Lemma'' below.  
 
 The ``Density Lemma'' is 
 

  \begin{lemma}\label{l.dense}  Any subset $\mathcal S\subset\mathcal T$ is a union of $\mathcal S_{\text{\rm heavy}}$ and $\mathcal S_{\text{\rm light}}$ for which 
  \begin{equation*}
 \dense {\mathcal S_{\text{\rm light}}}<\tfrac12\dense{ \mathcal S} ,
  \end{equation*}
 and the collection $\mathcal S_{\text{\rm heavy}}$  satisfies 
  \begin{equation} \label{e.densecount}
 \COUnt {\mathcal S_{\text{\rm heavy}}}\lesssim{}{\dense {\mathcal S}}^{-1}.
  \end{equation}
  \end{lemma}

  What is significant is that this relatively simple lemma admits a non-trivial variant  intimately
  linked to the tree structure and orthogonality.  We should refine the notion of a tree. Call a tree $\Tree$ with top $I_\Tree\times \omega_{\Tree}$
 a \tree\pm{} iff for each $s\in\Tree$, aside from the top, 
  $I_\Tree\times \omega_{\Tree}\cap I_s\times \omega_{s\pm}$ is not empty.  Any tree is a union of a 
  \tree+{} and a \tree-{}.  If $\Tree$ is a \tree+{}, 
  observe that the rectangles $\{I_s\times \omega_{s-}\mid s\in\Tree\}$ are disjoint.  And, by  the proof of 
  \Pcap.Qprop/, we see that 
   \begin{equation*}
   \Delta(\Tree)^2:=\sum_{s\in\Tree}\abs{\ipf }^2\lesssim{}\lVert f\rVert_2^2.
   \end{equation*}
  This motivates the  definition 
   \begin{equation} \label{e.size}
  \size {\mathcal S}:=\sup\{ \abs{I_\Tree}^{-1/2} \Delta(\Tree)\mid \Tree\subset\mathcal S,\ \text{$\Tree$ is a \tree+{}}\}.
   \end{equation}
  The ``Size Lemma'' is 
  
 
  \begin{lemma}\label{l.size}  Any subset $\mathcal S\subset\mathcal T$ is a union of $\mathcal S_{\text{\rm big}}$ and $\mathcal S_{\text{\rm small}}$ for which 
  \begin{equation*}
 \size {\mathcal S_{\text{\rm small}}}<\tfrac12\size {\mathcal S} ,
  \end{equation*}
 and the collection $\mathcal S_{\text{\rm big}}$ satisfies 
  \begin{equation} \label{e.sizecount} 
 \COUnt {\mathcal S_{\text{\rm big}}}\lesssim{}{\size {\mathcal S}}^{-2}.
  \end{equation}
  \end{lemma} 
 

 Our final Lemma relates trees, density and size. It is the ``Tree Lemma.'' 
 
  \begin{lemma}\label{l.tree}  For any tree $\Tree$ 
  \begin{equation}\label{e.tree}
 \sum_{s\in\Tree}\abs{\ipf \ip \phi_s,\ind E  .}\lesssim{}\abs{I_\Tree}\size \Tree \dense\Tree.
  \end{equation}
  \end{lemma}


 The final elements of the proof are organized as follows.  Certainly, $\dense {\mathcal T}<2$ for $\kappa$ sufficiently large.  We take 
 some finite subset $\mathcal S$ of $\mathcal T$, and so certainly $\size {\mathcal S}<\infty$.  If $\size {\mathcal S}<2$, we jump to the next stage of the proof. 
 Otherwise, we iteratively apply  \Lcap.size/ to obtain subcollections $\mathcal S_n\subset \mathcal S$, $n\ge0$,  
 for which 
  \begin{equation}\label{e.sizen}
 \size {\mathcal S_n}<2^{n},\quad  n>0,
  \end{equation}
 and $\mathcal S_n$  satisfies 
  \begin{equation} \label{e.count}
 \COUnt {\mathcal S_n} \lesssim{}2^{-2n}.
  \end{equation}
 We are left with a collection of tiles $\mathcal S'=\mathcal S-\bigcup_{n>0}\mathcal S_n$ which has both density and size at most $2$.  
 
 Now, both \Lcap.dense/ and \Lcap.size/ are set up for iterative application.  And we should apply them so that the estimates in 
 (\ref{e.densecount})  and (\ref{e.sizecount})  are of the same order. (This means that we should have density about the square of the size.)  As a consequence, we can achieve a decomposition of $\mathcal S$ into 
 collections $\mathcal S_n$, $n\in\mathbb Z$, which satisfy (\ref{e.sizen}) , (\ref{e.count})  and 
  \begin{equation}\label{e.densen}
 \dense {\mathcal S_n}<\min(2,2^{2n}).
  \end{equation}
 
 Use the estimates (\ref{e.tree}) ---(\ref{e.densen}) .  Write $\mathcal S_n$ as a union of trees $\Tree\in\widetilde{\Tree}_n$, this collection of trees 
satisfying the estimate of (\ref{e.count}) . We see that 
  \begin{align} \nonumber
 \sum_{s\in\mathcal S_n}\abs{\ipf \ip \phi_s,\ind E   .}={}&\sum_{\Tree\in\widetilde{\Tree}_n} \sum_{s\in\Tree}\abs{\ipf \ip \phi_s,{\mathbf 1_{E}  }.}
 \\  \label{e.ending}
 {}\lesssim{}&2^n\min(2,2^{2n})\sum_{\Tree\in\widetilde{\Tree}_n}\abs{I_\Tree}
 \\ \nonumber
 {}\lesssim{}&\min(2^{-n},2^n).
  \end{align}
 This is summable over $n\in\mathbb Z$ to an absolute constant, and so our proof (\ref{e.weaktypedual})  is complete, aside from the proofs of the three key lemmas.
 
 \subsection{Complements}

  \begin{remark}\label{r.weaktypeequiv} These two conditions are equivalent.  
  \begin{gather*}
 \sup_{\lambda>0} \lambda^{-2}\abs{\{ f>\lambda\}}\lesssim1,
 \\
 \int_E \abs f\; dx\lesssim\abs{E}^{1/2},\qquad \abs E<\infty.
  \end{gather*}
  \end{remark}

 \begin{remark}\label{r.dilate}  Let $\operatorname A$ be an operator for which $\dilate {2^k}^2 \operatorname A=\operatorname A
\dilate {2^k}^2$ for all $k\in\mathbb Z$.  
Suppose that 
there is an absolute constant $K$ so that for all functions $f\in L^2(\mathbb R)$ of norm one, 
 \begin{equation*}
\abs{ \{\operatorname  A f>1\}}\le{}K.
 \end{equation*}
Then for all $\lambda>0$, 
 \begin{equation*}
\abs{ \{\operatorname  Af>\lambda\}    }\lesssim{}\lambda^{-2}\lVert f\rVert_2^2.
 \end{equation*}
See \cite{MR37:6687}.
 \end{remark}

  \begin{remark}\label{r.sizemax}  For any $+$tree $\Tree$, 
  \begin{equation*}
 \sum_{s\in\Tree} \abs{\ip f,\varphi_s.}^2\lesssim \int \abs{f}^2 \trans {c(I_\Tree)}\dilate {I_\Tree}^\infty  \chi \; dx .
  \end{equation*}
  Moreover, one has the inequality 
  \begin{equation*}
  \abs{I_\Tree}^{-1}\sum_{s\in\Tree} \abs{\ip f,\varphi_s.}^2{}\lesssim{} \inf_{x\in\Tree}\operatorname M\abs{f}^2(x).
  \end{equation*}
 Here, $\operatorname M$ is the maximal function, 
  \begin{equation*}
\operatorname Mf(x)=\sup_{t>0}(2t)^{-1}\int_{-t}^t \abs{f(x-y)}\; dy.
  \end{equation*}
  \end{remark}

 
 \section{The Density Lemma}
 Set $\delta=\dense {\mathcal S}$. 
 Suppose for the moment that density had the simpler definition 
  \begin{equation*}
 \text{\sf dense}(s):= \frac{\abs{N^{-1}(\omega_{s})\cap I_s}}{\abs{I_s}}.
  \end{equation*}
 The collection $\mathcal S_{\text{\rm heavy}}$ is to be a union of trees.  So to select this collection, it suffices to select 
 the tops of the trees in this set. 
 
Select the tops of the trees, $\text{\sf Tops}$ as being those tiles $s$ with $\text{\sf dense}(s)$ exceeding 
 $\delta/2$, which are also maximal with respect to the partial order `$<$.'  The tree associated to such a tile $s\in \text{\sf Tops}$
 would just be 
 all those tiles in $\mathcal S$ which are less than $s$.  The   tiles in $\text{\sf Tops}$ are pairwise incomparable with respect to 
 the partial order `$<$,' and so are pairwise 
  disjoint rectangles in the time--frequency plane.  And so the sets ${N^{-1}(\omega_{s})\cap I_s}\subset E$ are pairwise disjoint, 
  and each has measure at  least $\frac\delta2\abs{I_s}$.  Hence the estimate below is immediate. 
  \begin{equation}\label{e.xtops}
 \sum_{s\in\text{\sf Tops}}\abs{I_s}\lesssim{}\delta^{-1}
  \end{equation}
 
 \smallskip

The Schwartz tails problem prevents us from using this very simple estimate to prove this lemma, but in the present context, 
the Schwarz tails  are a weak enemy at best.   Let ${\sf Tops}$ be those $s\in\mathcal S$ which have $\dense s>\delta/2$ and are maximal with
respect to `$<$.'  It suffices to show (\ref{e.xtops}) .
For an integer $k\ge0$, and  small constant $c$, let $\mathcal S_k$ be those $s\in{\sf Tops}$ for which 
 \begin{equation} \label{e.Bigg}
\abs{2^k I_s\cap N^{-1}(\omega_s)}\ge{}c2^{2k}\delta\abs{I_s}.
 \end{equation}
Every tile in ${\sf Tops}$ will be in some $\mathcal S_k$, with $c$ sufficiently small, and so it suffices to show that 
 \begin{equation}\label{e.dense2}
\sum_{s\in\mathcal S_k}\abs{I_s}\lesssim2^{-k}\delta^{-1}.
 \end{equation}
 
 Fix $k$.  Select from $\mathcal S_k$ a subset $\mathcal S'_k$ of tiles satisfying $\{2^kI_s\times \omega_s\mid s\in\mathcal S_k'\}$ are pairwise 
 disjoint, and if $s\in\mathcal S_k$ and $s'\in\mathcal S_k'$ are tiles such that $2^kI_s\times \omega_s $ and $2^kI_{s'}\times \omega_{s'}$ 
 intersect, then $\abs{I_s}\le\abs{I_{s'}}$.  It is clearly possible to select such a subset.  And since the tiles in $\mathcal S_k$ 
 are incomparable with respect to `$<$', we can use 
 (\ref{e.Bigg})  to estimate 
  \begin{align*}
 \sum_{s\in\mathcal S_k}\abs{I_s}\le{}&2^{k+1}\sum_{s'\in\mathcal S_k'}\abs{I_{s'}}
 \\{}\le{}&\tfrac2c2^{-k}\delta^{-1}.
  \end{align*}
 That is, we see that (\ref{e.dense2})  holds, completing our proof.

 \subsection{Complements} 
 
  \begin{remark}\label{r.bmo}  Let $\mathcal S$ be a set of tiles for which there is a constant $K$ so that for all dyadic intervals $J$, 
  \begin{equation*}
 \sum_{\substack{s\in\mathcal S\\ I_s\subset J}}\abs{I_s}\le{}K \abs J .
  \end{equation*}
 Then  for all $1\le{}p<\infty$, and intervals $J$, 
  \begin{equation*}
 \NOrm \sum_{\substack{s\in\mathcal S\\ I_s\subset J}}\ind {I_s}  .p. \lesssim{}K_p \abs {J}^{1/p} .
  \end{equation*}
 In fact, $K_p {}\lesssim{} p $. 
  \end{remark}

 \section{The Size Lemma}

Set $\sigma=\size {\mathcal S} $.   We will need to construct a collection of trees $\Tree\in{\widetilde\Tree}_{\text{large}}$, with 
$\mathcal S_{\text{large}}=\bigcup_{\Tree\in{\widetilde\Tree}_{\text{large}}}\Tree$, and 
 \begin{equation} \label{e.sizecnt}
\sum_{\Tree\in{\widetilde\Tree}_{\text{large}}}\abs{I_\Tree}\lesssim\sigma^{-2},
 \end{equation}
as required by (\ref{e.sizecount}) . 

The selection of trees $\Tree\in{\widetilde\Tree}_{\text{large}}$ will be done in conjunction with the 
construction of \tree+s $\Tree_+\in{\widetilde\Tree}_{\text{large}+}$.  This collection will play a critical role in the verification of (\ref{e.sizecnt}) . 

The construction is recursive in nature.  Initialize 
 \begin{equation*}
\mathcal S^{\text{stock}}:=\mathcal S,\qquad {\widetilde\Tree}_{\text{large}}:=\emptyset,\qquad {\widetilde\Tree}_{\text{large}+}:=\emptyset.
  \end{equation*}
 While $\size {{\mathcal S}^{\text{stock}} }>\sigma/2 $, select a \tree+{} $\Tree_+\subset\mathcal S^{\text{stock}}$ with 
  \begin{equation}\label{e.treeselect}
 \Delta(\Tree_+)>\frac \sigma 2 \abs{I_{\Tree_+}}.
  \end{equation}
 In addition, the top of the tree $I_{\Tree_+}\times \omega_{\Tree_+} $ should be maximal with respect to the partial order `$<$' 
 among all trees that satisfy (\ref{e.treeselect}) .  And $c(\omega_{\Tree_+})$ 
 should be minimal, in the order of $\mathbb R$. Then, take $\Tree$ to be the maximal tree (without reference to sign) in $\mathcal S^{\text{stock}}$ 
 with top $I_{\Tree_+}\times \omega_{\Tree_+} $.
 
 After this tree is chosen, update 
  \begin{gather*}
 \mathcal S^{\text{stock}}:=\mathcal S^{\text{stock}}-\Tree,\\
 {\widetilde\Tree}_{\text{large}}:={\widetilde\Tree}_{\text{large}}\cup \Tree,
 \qquad {\widetilde\Tree}_{\text{large}+}:={\widetilde\Tree}_{\text{large}+}\cup \Tree_+.
  \end{gather*}
 
 Once the while loop finishes, set $\mathcal S_{\text{small}}:=\mathcal S^{\text{stock}}$ and the recursive procedure stops. 
 
 \medskip 
 
 It remains to verify (\ref{e.sizecnt}) . This is a orthogonality statement, but one that is just weaker than true orthogonality. 
Note that a particular enemy is the is the situation in which $\ip \varphi_s,\varphi_{s'}.\not=0$.  When $\omega_s=\omega_{s'}$, 
this may happen, but as we saw in the proof of \Pcap.Qprop/, this case may be handled by direct methods. Thus we are primarily 
concerned with the case that e.g.~$\omega_{s-}\subset_{\not=}\omega_{s'-}$.

 A central part of this argument is a bit of geometry of the time--frequency plane that is encoded in 
 the construction of the \tree+s above.   Suppose there are two trees $\Tree\not=\Tree'\in{\widetilde\Tree}_{\text{large+}}$, and tiles $s\in\Tree$ and $s'\in\Tree'$ such that $\omega_{s-}\subset_{\not=}\omega_{s'-}$, 
 then, it is the case that $I_{s'}\cap I_{\Tree}=\emptyset$.  We refer to this property as `strong disjointness.' It is a condition that is 
 strictly stronger than just requiring that the sets in the time--frequency plane below are disjoint in $\Tree$.
  \begin{equation*}
 \bigcup_{s\in\Tree} I_s\times \omega_{s-},\qquad \Tree \in {\widetilde\Tree}_{\text{large}+}
  \end{equation*}
 
 To see that strong disjointness  holds,  observe that $\omega_\Tree\subset \omega_s\subset_{\not=} \omega_{s'-}$.
 Thus $\omega_{\Tree'}$ lies above $\omega_{\Tree}$.  That is, in our recursive procedure, $\Tree$ was constructed first.  
 If it were the case that 
 $I_{s'}\cap I_{\Tree}\not=\emptyset$, observe that one interval would have to be contained in the other. 
 But tiles have area one, thus, it 
 must be the case that $I_{s'}\subset I_\Tree$.  That means that $s'$ would have been in the tree (the one without sign) 
 that was removed from 
 $\mathcal S^{\text{stock}}$ before $\Tree'$ was constructed.  This is a contradiction which proves strong disjointness. 
See \f.strongly-disjoint/.
 
\begin{figure}
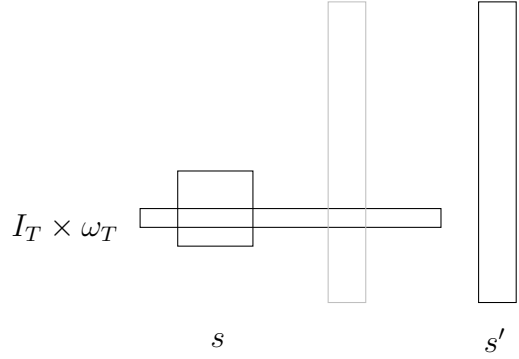
 
\begin{center}
\begin{pgfpicture} {0cm}{0cm}{4cm}{4.7cm}

\color{lightgray}

\color{black} 
\pgfrect[stroke]{\pgfxy(1,2)}{\pgfxy(4,.25)} 
\pgfrect[stroke] {\pgfxy(1.5,1.75)}{\pgfxy(1,1)} 
\pgfrect[stroke]{\pgfxy(5.5,1,)}{\pgfxy(.5,4)} 

\pgfputat{\pgfxy(0,2.17)}{\pgfbox[center,top]{$I_{T}\times\omega_{T}$} }
	
\pgfputat{\pgfxy(2.1,.5)}{\pgfbox[center,center]{$s$ } }

\pgfputat{\pgfxy(6,.5)}{\pgfbox[right,center]{$s'$ } }

\color{lightgray}  
\pgfrect[stroke]{\pgfxy(3.5,1,)}{\pgfxy(.5,4)}

\end{pgfpicture}
\end{center}
\label{f.strongly-disjoint} 
	\caption{The proof of strongly disjoint trees. Note that the gray tile could be in the tree that was removed after the 
	selection of the $+$--tree with the top indicated above. } 
\end{figure} 
 \smallskip 
 
 We use this strong disjointness condition, and the selection criteria (\ref{e.treeselect}) , to prove the bound (\ref{e.sizecnt}) .  The method of proof 
 is closely related to the so called $\operatorname T \operatorname T^*$ method.    Set $\mathcal S'=\bigcup_{\Tree\in {\widetilde\Tree}_{\text{large}+}}\Tree$, and
  \begin{equation*}
 F:=\sum_{s\in\mathcal S'}\ip f,\varphi_s. \varphi_s .
  \end{equation*} 
 
The operator $f\mapsto \ip f,\varphi_s. \varphi_s$  is self--adjoint, so that 
  \begin{align*}
 \sigma^2\sum_{\Tree\in {\widetilde\Tree}_{\text{large}+}}\abs{I_\Tree}={}& 
 \ip f, F. 
 \\
 {}\le{}& \norm f.2. \norm F .2.
  \end{align*}
 And so, we should show that 
  \begin{equation} \label{e.Stree}
 \lVert F\rVert_2^2\lesssim{} 
 \sigma^2\sum_{\Tree\in {\widetilde\Tree}_{\text{large}+}}\abs{I_\Tree}.
  \end{equation}
 This will complete the proof. 
 
 \smallskip 
 This last inequality is seen by expanding the square on the left hand side.  In particular, the left hand side of (\ref{e.Stree})  is 
 at most the sum of the two terms
  \begin{gather} \label{e.equals} 
 \sum_{\substack{s,s'\in\mathcal S'\\ \omega_{s}=\omega_{s'}} }\ip f,\varphi_s. \ip \varphi_s ,\varphi_{s'}. \ip \varphi_{s'},f. 
 \\
 \label{e.strictsubset}
 2\sum_{\substack{s,s'\in\mathcal S'\\ \omega_{s}\subset_{\not=}\omega_{s'}}}\abs{ \ip f,\varphi_s. \ip \varphi_s ,\varphi_{s'}. \ip \varphi_{s'},f.}
  \end{gather}
 
 For the term  (\ref{e.equals}) , we have the obvious estimate on the inner product 
  \begin{equation*}
 \abs{ \ip \varphi_s ,\varphi_{s'}.}\lesssim{} \Bigl(1+ \frac{\operatorname{dist}(I_s,I_{s'})}{\abs{I_s}}\Bigr)^{-4}.
  \end{equation*}
  (Compare to (\ref{e.compare}) .)  Thus, by Cauchy--Schwarz, 
  \begin{equation*}
(\ref{e.equals}) \lesssim{} \sum_{s\in\mathcal S'} \abs{\ip f,\varphi_s.}^2\lesssim{}\sigma^2\sum_{\Tree\in {\widetilde\Tree}_{\text{large}+}}\abs{I_\Tree}. 
 \end{equation*}

For the term (\ref{e.strictsubset}) ,  we need only show that for each tree $\Tree$, 
 \begin{equation}  \label{e.laststep}
\sum_{s\in\Tree } { \sum_{\substack{s'\in\mathcal S'\\ \omega_{s}\subset_{\not=}\omega_{s'}}}
\abs{\ipf \ip \varphi_s,\varphi_{s'}. \ip \varphi_{s'},f.}}^2\lesssim{}\sigma^2\abs{I_\Tree}.
 \end{equation}
Here,  $\mathcal S(s):=\{s'\in\mathcal S'-\Tree\mid \omega_{s-}\subset_{\not=}\omega_{s'-}\}$. 
The implied constant should be independent of the tree $\Tree$.

Now, the strong disjointness condition enters in two ways.  For $s\in\Tree$, and $s'\in\mathcal S(s)$, it is the case that 
$I_{s'}\cap I_{\Tree}=\emptyset$.  But furthermore, for $s',s''\in\mathcal S(s)$,
we have e.g.~$\omega_{s-}\subset \omega_{s'-}\subset \omega_{s
''-}$, 
so that $I_{s'}\cap I_{s''}$ is also empty. 

At this point,  rather clumsy   estimates of (\ref{e.laststep})  
are in fact optimal.  The definition of size gives us the bound 
 \begin{equation*} 
\abs{\ip \varphi_{s'},f.}\lesssim{} \sqrt{\abs{I_{s'}}}\sigma
 \end{equation*}
And, since $\omega_{s}\subset \omega_{s'}$, we have $\abs{I_s}\ge\abs{I_{s'}}$, and $I_s$, and $I_{s'}$ are,
in the typical situation, 
far apart.  An  estimation left to the reader gives 
 \begin{equation}  \label{e.lefttoreader}
\abs{\ip \varphi_s,\varphi_{s'}.}\lesssim{}{\sqrt{\abs{I_{s'}}\abs{I_{s}}}} \chi_{I_s}(c(I_{s'}))
 \end{equation}
Thus, we bound the left side of (\ref{e.laststep})  by 
 \begin{align} \notag
\sum_{s\in\Tree }  \sum_{s'\in\mathcal S(s)} \abs{ \ipf \ip \varphi_s,\varphi_{s'}. \ip \varphi_{s'},f.}\lesssim{}&
\sigma^2 \sum_{s\in\Tree }  {\abs{I_{s'}}} \abs{I_s}  \chi_{I_s}(c(I_{s'})) 
\\{}\lesssim{}& \notag
\sigma^2 \sum_{s\in\Tree} \int_{(I_\Tree)^c} \abs{I_s}\chi_{I_s}(x) \;dx 
\\{}\lesssim{}& \sigma^2 \abs{I_\Tree}   \label{e.easytoverify}
  \end{align}
as is easy to verify. This completes the proof of (\ref{e.laststep}) , and so finishes the proof of \Lcap.size/.

\subsection{Complements}

 \begin{remark}\label{r.easytoverify}
Concerning the inequality (\ref{e.easytoverify}) ,  for any tree $\Tree $, we have 
 \begin{equation*}
\sum _{s\in \Tree} \int _{I_\Tree^c} \abs{I_s}\chi_{I_s}(x) \; dx {}\lesssim{}\abs {I_\Tree}.
 \end{equation*}
 \end{remark}

 \begin{remark}\label{r.tree2}    
 Let $\Tree$ be a $+$tree and set 
 \begin{equation*}
F_\Tree=\sum_{s\in\Tree}\ip f,\varphi_s. \varphi_s 
 \end{equation*}
Then, the inequality below is true. 
 \begin{align*}
\norm F_\Tree .2.\simeq{}&\Bigl[ \sum_{s\in\Tree}\abs{\ip f,\varphi_s.}^2 \Bigr]^{1/2}
\\{}\lesssim{}& \size \Tree \abs{I_\Tree}^{1/2}.
 \end{align*}
 \end{remark}

 \begin{remark}\label{r.bmosize}   With the notation above, assume that $0\in\omega_\Tree $.  Then, 
 \begin{align*}
\size \Tree {}\simeq{} & \sup_J \bigl[ \abs{J}^{-1}\sum_{\substack{s\in\Tree\\ I_s\subset J}} \abs{\ip f,\varphi_s .}^2 
	\bigr]^{1/2} \\
{}\simeq{}& \sup_J \Bigl[ \abs{J}^{-1}\int_J \ABs{ F_\Tree - \abs{J}^{-1}\int_J F_\Tree }^2 dx\Bigr]^{1/2},
 \end{align*}
where the supremum is over all intervals $J$.  The last quantity is the $BMO$ norm of $F_\Tree$. 
 \end{remark}

 \begin{remark}\label{r.weaktype}  
It is an important heuristic that for a collection $\mathcal S$ of pairwise incomparable tiles, the 
functions $\{\varphi_s\mid s\in\mathcal S\}$ are nearly orthogonal.   The heuristic permits a quantification in 
terms of the following weak type inequality. 
Let $\mathcal S$ be a collection of tiles that are pairwise incomparable with respect to 
`$<$.'  Then for all $f\in L^2$ and all $\lambda>0$, 
 \begin{equation*}
\sum_{s\in\mathcal S_\lambda}\abs {I_s} {}\lesssim{} \lambda^{-2}\lVert f\rVert_2^2,
 \end{equation*}
where $\mathcal S_\lambda=\{s\in\mathcal S\mid \abs{\ip f,\varphi_s.}>\lambda\sqrt{\abs{I_s}}\}$. 
Note that this in an inequality about the boundedness of a sublinear operator from $L^2(\mathbb R)$ to $L^{2,\infty}(\mathbb R\times\mathcal S)$.
In the latter space, one uses counting measure on $\mathcal S$.
 \end{remark}

  \begin{remark}\label{r.weaktypetree}
 Another important heuristic is that the notion of ``strong disjointness'' for trees is as ``pairwise incomparable''
 is for tiles. 
 Let $\widetilde\Tree$ be a collection of strongly disjoint trees.  Show that for all $f\in L^2$ and 
 all $\lambda>0$, 
  \begin{equation*}
 \sum_{\Tree\in\widetilde\Tree_\lambda}\abs{I_\Tree} {}\lesssim{} \lambda^{-2}\lVert f\rVert_2^2,
 \end{equation*}
where $\widetilde\Tree_\lambda=\{ \Tree \in \widetilde\Tree\mid \Delta(\Tree) \ge\lambda\abs{I_\Tree}\}$. 
 \end{remark} 

 \begin{remark}\label{r.above2}  
 Let $\mathcal S$ be a collection of tiles that are pairwise incomparable with respect to 
`$<$.'  Show that for all $2<p<\infty$, 
 \begin{equation*}
\Bigl[ \sum_{s\in\mathcal S} \abs{\ip f ,\varphi_s .}^p\Bigr]^{1/p}
{}\lesssim{} \norm f.p. . 
 \end{equation*}
Notice that the form of this estimate at $p={} \infty $ is obvious. 
 \end{remark}

 \begin{remark}\label{r.ABove2}  
 Let $\widetilde\Tree$ be a collection of strongly disjoint trees. 
Then for all $2<p<\infty$, 
 \begin{equation*}
  \sum_{\Tree\in \widetilde\Tree} {\Delta(\Tree)^p} 
{}\lesssim{} \norm f.p. . 
 \end{equation*}
 \end{remark}

 \begin{remark}\label{r.rubio}   
The $L^p$ estimates of the previous two complements can in some instances be improved.  
For each integer $k$, 
 \begin{equation*}
\NOrm \Bigl[ \sum_{\substack{s\in\mathcal T \\ \abs{I_s}=2^k}} \frac{ \abs{\ip f,\varphi_s.}^2} {\abs{I_s}} \ind {I_s} \Bigr]^{1/2} .p.
{}\lesssim{} \norm f.p.\qquad 2<p<\infty. 
 \end{equation*}
This can be seen by  showing that 
 \begin{equation*}
\sum_{\substack{s\in\mathcal T \\ \abs{I_s}=2^k}} \frac{ \abs{\ip f,\varphi_s.}^2} {\abs{I_s}} \ind {I_s}
{}\lesssim{}(\dilate {2^k}^1 \chi)* \abs{f}^2.
 \end{equation*}
 \end{remark}

\section{The Tree Lemma}

We begin with some remarks about the maximal function, and a particular form of the same that we shall use 
at a critical point of this proof.  Consider the maximal function 
 \begin{equation*}
\operatorname Mf=\sup_{I\in\mathcal D}\ind I \abs{\ip f,\chi_I.}.
 \end{equation*}
It is well known that this is bounded on $L^2$.  A proof follows.  Consider a linearized version of the supremum.  
To each $I\in\mathcal D$, associate a set $E(I)\subset I$, and require that the sets $\{E(I)\mid I\in\mathcal D\}$ be pairwise disjoint. 
(Thus, for fixed $f$, $E(I)$ is that subset of $I$ on which the supremum above is equal to $\abs{\ip f,\chi_I.}$.) Define 
 \begin{equation*} 
\operatorname A f=\sum_{I\in\mathcal D}\ind {E(I)} \ip f,\chi_I.
 \end{equation*}
We show that $\norm \operatorname A .2.$ is bounded by a constant, independent of the choice of the sets $E(I)$.  

The method is that of $\operatorname T \operatorname T^*$.  Note that for positive $f$
 \begin{align*}
\operatorname A\operatorname A^*f\le{}&2\sum_{I\in\mathcal D}\sum_{\abs J\le\abs I} \ind {E(I)} \ip \chi_I,\chi_J.\ip \ind {E(J)} ,f.
\\{}\lesssim{}& \sum_{I\in\mathcal D}\ind {E(I)} \ip f,\chi_I.
 \end{align*}
It follows that 
 \begin{align*}
\lVert\operatorname A^*\rVert_2^2={} &  \sup_{\norm f.2.=1} \ip \operatorname A^* f ,\operatorname A^*f .
\\{}={}&\sup_{\norm f.2.=1}\ip f,\operatorname A\operatorname A^* f. 
\\{}\lesssim{}&{}\norm \operatorname A.2.
 \end{align*}
and so $\norm \operatorname A.2.\lesssim{}1$, as claimed.

We shall have recourse to not only this, but a particular refinement.
Let $\mathcal J$ be a partition of $\mathbb R$ 
into dyadic intervals.  To each $J\in\mathcal J$, associate a subset $G(J)\subset J$, with $\abs{G(J)}\le{}\delta\abs J$, 
where $0<\delta<1$ is fixed. 
Consider 
 \begin{equation}  \label{e.maxWithDense}
\operatorname M_\delta f:=\sum_{J\in\mathcal J}\ind {G(J)} \sup_{I\supset J}  \abs{\ip f,\chi_I.}
 \end{equation}
Then $\norm \operatorname M_\delta .2.\lesssim{}\sqrt\delta$.  The proof is  
 \begin{align*}
\int \abs{\operatorname M_\delta f}^2\; dx={}& \sum_{J\in\mathcal J}\abs{G(J)}\sup_{I\supset J}  \abs{\ip f,\chi_I.}
\\{}\le{}&\delta \sum_{J\in\mathcal J}\abs{J}\sup_{I\supset J}  \abs{\ip f,\chi_I.}
\\{}\le{}&\delta \int \abs{\operatorname M f}^2\; dx
 \\{}\lesssim{}&\delta\lVert f\rVert_2^2.
  \end{align*}
 
 \bigskip 
 We begin the main line of the argument. Let $\delta=\dense \Tree $, and $\sigma=\size \Tree $. Make a choice of signs $\varepsilon_s\in\{\pm1\}$ such that 
  \begin{equation*}
 \sum_{s\in\Tree}\abs{\ip f,\varphi_s. \ip \phi_s,\ind E .}={} 
 	\int_E{\sum_{s\in\Tree} \varepsilon_s \ip f,\varphi_s. \phi_s }\; dx.
 \end{equation*}
 By the ``Schwartz tails,'' the integral above is supported on the whole real line.  
 Let $\mathcal J$ be a partition of $\mathbb R$ consisting of the maximal dyadic intervals $J$ such that $3J$ does not contain any $I_s$ for $s\in\Tree$.  
 It is helpful to observe that for such $J$ if $\abs{J}\le{}\abs{I_\Tree}$, then $J\subset 3 I_\Tree$. And if $\abs J\ge\abs {I_\Tree}$,
 then 
 $\text{dist}(J,I_\Tree)\gtrsim\abs J$. 
 The integral above is at most the sum over $J\in\mathcal J$ of the two terms below.
  \begin{gather} \label{e.firstsum}
 \sum_{\substack{s\in\Tree\\ \abs{I_s}\le\abs J}}\abs{\ip f,\varphi_s. }\int_{J\cap E}\abs{\phi_s}\; dx
 \\
 \label{e.secondsum}
 \int_{J\cap E} \ABs{ \sum_{\substack{s\in\Tree\\ \abs{I_s}>\abs J}}\varepsilon_s\ip f,\varphi_s.  \phi_s }\; dx 
  \end{gather}
 Notice that for the second sum to be non--zero, we must have $J\subset 3I_\Tree$.

 The first term (\ref{e.firstsum})  is controlled by an appeal to the ``Schwartz tails.''  Fix an integer $n\ge0$, and only consider those $s\in\Tree$ for which $\abs{I_s}=2^{-n}\abs J$. Now, the distance of $I_s$ to $J$ is at least ${}\gtrsim\abs J$.  An
d, 
 \begin{equation*}
\abs{\ip f,\varphi_s.} \int_{J\cap E}\abs{\phi_s}\; dx\le{}\sigma\delta (\abs{I_s}^{-1}\text{dist}(I_s,J))^{-10} \abs{I_s}.
 \end{equation*}
The $I_s\subset I_\Tree$, so that summing this over $\abs{I_s}=2^{-n}\abs J$ will give us 
 \begin{equation*}
\sigma\delta 2^{-n}\min\bigl(\abs J, \abs {I_\Tree} (\text{dist}(J,I_\Tree)\abs{I_\Tree}^{-1})^{-10}\bigr).
 \end{equation*}
This is summed over $n\ge0$ and $J\in\mathcal J$ to bound (\ref{e.firstsum})  by ${}\lesssim{}\sigma\delta\abs{I_\Tree}$, as required.

\medskip

Critical to the control of (\ref{e.secondsum})  is the following observation. Let 
 \begin{equation} \label{e.GJ} 
G(J)=J\cap\bigcup_{\substack{s\in\Tree\\ \abs{I_s}\ge\abs J}} N^{-1}(\omega_{s+}).
 \end{equation}
Then $\abs{G(J)}\lesssim{}\delta\abs J$. To see this, let $J'$ be the next larger dyadic interval that contains 
$J$.  Then $3J'$ must contain some $I_{s'}$, for $s'\in\Tree$.  Let $s''$ be that tile with $I_{s'}\subset I_{s''}$, 
$\abs{I_{s''}}=\abs J$, and $\omega_\Tree\subset \omega_{s''}$.  Then, $s'<s''$, and by the definition of density, 
 \begin{equation*}
\int_{E\cap N^{-1}(\omega_{s''})} \chi_{I_{s''}} \; dx\le\delta
 \end{equation*}
But, for each $s$ as in (\ref{e.GJ}) , we have $\omega_s\subset \omega_{s''}$, so that
 $G(J)\subset N^{-1}(\omega_{s''})$. Our claim follows.

 Suppose that $\Tree$ is a $-$tree.  That means that the tiles $\{I_s\times \omega_{s+}\mid s\in\Tree\}$ are disjoint.  
 We use an estimation absent of any cancellation effects.  
 Then, the bound for (\ref{e.secondsum})  is no more than 
  \begin{equation*}
 \abs{G(J)} \NOrm \sum_{\substack{s\in\Tree\\ \abs{I_s}\ge\abs J}} \abs{\langle f,\varphi_s\rangle \phi_s} .\infty.\lesssim{}\delta\sigma\abs{J}.
  \end{equation*}
 This is summed over $J\subset 3I_\Tree$ to get the desired bound.

 Suppose that $\Tree$ is a $+$tree.  (This is the interesting case.)
 Then, the tiles $\{I_s\times \omega_{s-}\mid s\in\Tree\}$ are pairwise disjoint, and we set 
  \begin{equation*}
 F=\modulate {-c(\omega_\Tree)}\sum_{s\in\Tree}\varepsilon_s\ip f,\varphi_s. \varphi_s 
  \end{equation*}
 Here it is useful to us that we only use the ``smooth'' functions $\varphi_s$ in the definition of this function. 
 Note that $\norm F.2.\lesssim{}\sigma\sqrt{\abs{I_\Tree}}$, which is a consequence of the definition of size and \Pcap.Qprop/.  
 Set $\tau(x)=\sup\{\abs{I_s}\mid s\in \Tree,\ N(x)\in\omega_{s+}\}$, and observe that 
 for each $J$, and $x\in J$, 
  \begin{equation*}
 \sum_{\substack{s\in\Tree\\ \abs{I_s}\ge\abs J}}\varepsilon_s\ip f,\varphi_s. \phi_s(x)={}
 	\sum_{\substack{s\in\Tree\\ \tau(x)\ge\abs{I_s}\ge\abs J}}\varepsilon_s\ip f,\varphi_s. \varphi_s(x)
  \end{equation*}
 This is so since all of the intervals $\omega_{s+}$ must contain $\omega_\Tree$, and  if $N(x)\in\omega_{s+}$, then it must also 
 be in every other $\omega_{s'+}$ that is larger.  What is significant here is that on the right we have a truncation of the sum 
 that defines $F$. 
 
 This last sum can be dominated by a maximal function.  For any $\tau>0$ and $J\in\mathcal J$, let 
  \begin{equation*}
 F_{\tau,J}=\modulate {-c(\omega_\Tree)}\sum_{\substack{s\in\Tree\\ \tau\ge\abs{I_s}\ge\abs J}}\varepsilon_s\ip f,\varphi_s. \varphi_s 
 \end{equation*}
This function has Fourier support in the interval $[-\frac78\abs{J}^{-1},-\frac14\tau^{-1}]$.  In particular, recalling how we defined 
$\varphi$, we can choose $\frac1{16}<a,b<\frac14$ so that 
 \begin{equation*}
 F_{\tau,J}=(\dilate {a\abs J}^1 \varphi {}-\dilate {b\tau} ^1 \varphi) * F
  \end{equation*}
 We conclude that for $x\in J$, 
  \begin{equation*}
 \abs{ F_{\tau(x),J}(x)}\lesssim{} \operatorname M_\delta F (x),
  \end{equation*} 
 where $\operatorname M_\delta$ is defined as in (\ref{e.maxWithDense}) .  
 
 The conclusion of this proof is now at hand. We have 
  \begin{align*}
 \sum_{\substack{J\in\mathcal J\\ \abs{J}\le3abs{I_\Tree}}} 
 \int_{G(J)}\abs{F_{\tau(x),J}}\; dx {}\lesssim{}& \int_{\bigcup_{\abs{J}\le3\abs{I_\Tree}} G(J)} \operatorname M_\delta F\; dx
 \\{}\lesssim{}& \ABs{\bigcup_{\abs{J}\le3\abs{I_\Tree}} G(J)}^{1/2} \norm\operatorname M_\delta F.2.
 \\{}\lesssim{}& \delta\sqrt{\abs{I_\Tree}} \norm  F.2.
 \\{}\lesssim{}& \sigma\delta {\abs{I_\Tree}} 
 \end{align*}

\subsection{Complements} 

 \begin{remark}\label{r.simpletree}  The estimate below is somewhat cruder than the one just obtained, and therefore 
easier to obtain. For all trees $\Tree$, 
 \begin{align*}
\NOrm \sum_{s\in\Tree }\langle f,\varphi_s\rangle \phi_s .2.\lesssim{}& \Bigl[ \sum_{s\in\Tree} \abs{\ip f,\varphi_s.}^2\Bigr]^{1/2}
\\{}\lesssim{}& \size \Tree \abs{ I_\Tree}^{1/2}.
 \end{align*}
 \end{remark} 

 \begin{remark}\label{r.max}  The  maximal function $\operatorname M_\delta$ in (\ref{e.maxWithDense})  admits the bounds 
 \begin{equation*}
\norm \operatorname M_\delta.p.\lesssim{} \delta^{1/p},\qquad 1<p<\infty.
 \end{equation*}
This depends upon the fact that the maximal function itself maps $L^p$ into itself, for $1<p<\infty$. 
 \end{remark}

 \begin{remark}\label{r.treep}  For any tree $\Tree$, 
 \begin{equation*}
\NOrm \sum_{s\in\Tree}\langle g,\varphi_s\rangle \phi_s .p.\lesssim{} \delta^{1/p}\norm g.p.,\qquad 1<p<\infty.
 \end{equation*}
 \end{remark} 

 \begin{remark}\label{r.LP}  For a $+$tree $\Tree$, 
 \begin{equation*}
\NOrm \sum_{s\in\Tree} \langle f,\varphi_s\rangle\varphi _s.p.\lesssim{}\NOrm \Bigr[\sum_{s\in\Tree} \frac{\abs{\ip f,\varphi_s.}^2}{\abs{I_s}} \ind {I_s} \Bigr]^{1/2} .p.,\qquad 1<p<\infty. 
 \end{equation*}
Conclude that 
 \begin{equation*}
\NOrm \sum_{s\in\Tree} \langle f,\varphi_s\rangle\varphi_s .p.\lesssim{} \size \Tree \abs{I_\Tree}^{1/p}, 
\qquad 1<p<\infty. 
 \end{equation*}
 \end{remark}


\section{Carleson's Theorem on $L^p$, $1<p\not=2<\infty$ }
 
 We outline a proof that the Carleson maximal operator maps $L^p$ into itself for all $1<p<\infty$. 
 The key point is that we should obtain a distributional estimate for the model operator.

  \begin{proposition}\label{p.model_pnot2}  
 For $1<p<\infty$, there is an absolute constant $K_p$ so that for all sets $E\subset\mathbb R$ of finite 
 measure and measurable functions $N$, we have 
  \begin{equation} \label{e.distributional}
 \abs{ \{ \abs{ \mathcal C_N \ind E } >\lambda\}} {}\le{} K_p^p\lambda^{-p}\abs{E}.
  \end{equation} 
  \end{proposition} 

Interpolation provides the $L^p$ inequalities.  We shall in fact prove that for all sets $E$  and $F$, there 
is a set $F'\subset F$ of measure $\abs{F'}\ge\frac12\abs F$,  
 \begin{equation} \label{e.EFbest}
\abs{\ip \mathcal C \ind E ,\ind {F'} .}\lesssim{} \min( \abs E ,\abs F )( 1+ \abs{ \log \tfrac{\abs E}{\abs F}} )
 \end{equation}
It is a routine matter to see that this estimate implies that 
 \begin{equation} \label{e.best}
\abs{ \{ \abs{\mathcal C \ind E }>\lambda\} }\lesssim{} \abs E 
\begin{cases} 
	\lambda\abs{ \log \lambda} & \text{if $ 0<\lambda<1/2$  }\cr
	e^{-c\lambda} 	& \text{otherwise} \cr
\end{cases}
 \end{equation}
Here, $c$ is an absolute constant.  This distributional inequality is in fact the best that is known about the 
Carleson operator.  See Sj\"olin \cite{MR39:3222}   
for the Walsh case, and \cite{MR49:998} for the Fourier case. 
A more recent publication proving the same point is  Arias de Reyna \cite{MR1906800}.  Both authors present a proof 
 along the lines of   Carleson.  We follow the weak type inequality  approach of Muscalu, 
Tao, and Thiele 
\cite{MR2003b:42017}.  The relevance of this approach to the Carleson theorem was demonstrated by Grafakos, Tao, and Terwilleger 
\cite{MR2031458}.

We shall find it necessary to appeal to some deeper properties of the Calder\'on Zygmund theory, and in particular 
a weak $L^1$ bound for the maximal function, but also the  bound  in (\ref{e.cz})  below.  

In proving (\ref{e.EFbest}) , we can rely  upon invariance under dilations, up to 
a change in the measurable $N(x)$, to assume that $1/2<\abs E \le1$.  As we already know the weak $L^2$ estimate, 
(\ref{e.EFbest})  is obvious for $\frac13<\abs F<3$.  The argument then splits into two 
cases, that of $\abs{F}<\frac13$ or $\abs F\ge3$.  

Note that our measurable function $N(x)$ is defined on the set $F$.  It is clear that our Density Lemma, \l.dense/, continues 
to hold in this context, with the change that the measure of $F$ should be added to the right hand side of (\ref{e.densecount}) .

 \subsection{The case of  $\abs F<\frac13$} 
 In this case, we will take $F'=F$. 
 Recall that $\mathcal T$ denotes the set of all tiles. 
  Clearly, $\size {\mathcal T} {}\lesssim1$.  We repeat the argument of 
 (\ref{e.sizen}) ---(\ref{e.ending}) .   Here, we should keep in mind that we want to balance out the estimate for the $\COUnt \cdot $ function, 
 and that we have a better upper bound on the count function coming from the Density Lemma.  
 Thus, $\mathcal T$ is a union of  collections $\mathcal S_n$, for $n\ge0$, so that 
  \begin{gather}
 \label{e.Dense} 
 \dense {\mathcal S_n} {}\lesssim2^{-2n},
 \\
 \label{e.Size} 
 \size {\mathcal S_n} {}\lesssim{}\min(1,2^{-n}\abs{F}^{-1/2}),
\\
 \label{e.Count} 
 \COUnt {\mathcal S_n} {}\lesssim2^{2n}\abs F .
  \end{gather}
 Then by the calculation of (\ref{e.ending}) , we have 
  \begin{align*}
 \sum_{s\in\mathcal S_n}\abs{\ip \ind E ,\varphi_s . \ip \phi_s, \ind {F'} .}\lesssim{} & 
 	\dense  {\mathcal S_n}  \size  {\mathcal S_n} 2^{2n}\abs F 
\\{}\lesssim{}& \min( \abs F ,\abs{F}^{1/2} 2^{-n}).
 \end{align*}
The sum of this terms over $n\ge0$ is no more than 
 \begin{equation*}
{}\lesssim\abs F \cdot\abs{ \log \abs F }
 \end{equation*}
which is as required. 

\subsection{The case of $\abs F\ge3$}
  
This case corresponds to the analysis of the Carleson operator on $L^p$ for $1<p<2$.  We shall have need of a 
more delicate weak type inequality below to complete this proof. 
To define the set $F'$,  let 
   \begin{equation*}
  \Omega=\{ \operatorname M\ind E > C_1 \abs{F}^{-1}\}.
   \end{equation*}
By the weak $L^1$ inequality for the maximal function, for an absolute  choice of $C_1$, we have 
  $\abs{\Omega}<\frac12\abs F$.  And we take $F'=F\cap\Omega^c$. The inner product in (\ref{e.EFbest})  is less than the sum of  
   \begin{gather}  \label{e.inside} 
  \sum_{\substack{s\in\mathcal T\\ I_s\subset\Omega}} \abs{\ip f,\varphi_s.\ip \phi_s,\ind {F'} .} 
\\    \sum_{\substack{s\in\mathcal T\\ I_s\not\subset\Omega}} \abs{\ip f,\varphi_s.\ip \phi_s,\ind {F'} .} 
\label{e.outside}
 \end{gather}
These sums are handled separately.

\bigskip

For (\ref{e.inside}) ,  observe that $ \varphi_s $ is essentially supported inside of $\Omega$ while $\phi_s$ is essentially not supported there. 
Thus, we should rely upon Schwartz tails to handle this term.
Let $J\subset \Omega$ be an interval such that $2^kJ\subset \Omega$ but $2^{k+1}J\not\subset\Omega$.  We observe two inequalities for such an interval, which 
are stated using the function $\chi_J$, as defined in (\ref{e.chiI}) .  The first 
is that 
 \begin{equation*}
\int_{F'}\chi_J\; dx\le{}\int_{(2^kJ)^c}\chi_J \; dx\lesssim{}2^{-(\kappa-1) k}
  \end{equation*}
Here, $\kappa$ is a large constant in the definition of $\chi$.   Also, we have 
 \begin{align*} 
\int_E \chi_J \; dx\lesssim{}& 2^{k}\int_E \chi_{2^{k+1}J}\; dx 
 \\{}\lesssim{}&  2^{k}\inf_{x\in 2^{k+1}J}\operatorname M\ind E (x)
   \\{}\lesssim{}& 2^{k} \abs{F}^{-1}.
   \end{align*}
  The last line follows  as some point in $2^{k+1}J$ must be in $\Omega$.   
  
  Observe that among all tiles $s$ with $I_s=J$, there is exactly one tile $s$ with 
  $N(x)\in\omega_{s+}$.  Hence
   \begin{align*}
 \sum_{\substack{s\in\mathcal T\\I_s=J}}\abs{\ip f,\varphi_s.\ip \phi_s,\ind {F'} .}
 {}\lesssim{}& \abs J\sum_{\substack{s\in\mathcal T\\I_s=J}} \ip \ind E , \chi_J. \ip \ind {F'} ,\chi_J .
 \\{}\lesssim{}& 2^{-(\kappa -2)k}\abs{F}^{-1}\abs J .
  \end{align*}
 Recall that $k$ is associated to how deeply $J$ is embedded in $\Omega$, and that $\Omega$ has measure 
 at most ${}\lesssim\abs{F}$.  Hence the right hand side above 
 can be summed over $J\subset\Omega$ to see that 
  \begin{equation*}
 \text{(\ref{e.inside}) }\lesssim{}1,
  \end{equation*}
which is better than  desired.

 \bigskip
 We turn to the second estimate. 
 Set $\mathcal T_{\text{out}}:=\{s\in\mathcal T\mid I_s\not\subset\Omega\}$, which is the collection of 
 tiles summed over in (\ref{e.outside}) . 
 The essential aspect of the definition of $\Omega$  is this Lemma.
 \begin{lemma}\label{l.sizesmall} 
  \begin{equation*}
\size {\mathcal T_{\text{\rm out}}}  {}\lesssim{}\abs{F}^{-1}.
 \end{equation*}
 \end{lemma}

Assuming the Lemma, we turn to the line of argument (\ref{e.sizen}) ---(\ref{e.ending}) .  The collection $\mathcal T_{\text{out}}$ can be decomposed into 
collections $\mathcal S_n$, for $n\ge0$, for which (\ref{e.Dense})  and (\ref{e.Count})  holds, and in addition 
 \begin{equation*}
\size {\mathcal S_n} {}\lesssim{}\min( \abs{F}^{-1}, \abs{F}^{-1/2}2^{-n}).
 \end{equation*}
Then by the calculation of (\ref{e.ending}) , we have 
  \begin{equation*}
 \sum_{s\in\mathcal S_n}\abs{\ip \ind E ,\varphi_s . \ip \phi_s, \ind {F'} .}\lesssim{} \min( 1 ,\abs{F}^{1/2} 2^{-n}),
 \end{equation*}
 making the sum over $n\ge0$ no more than $\log \abs F $, as required. This completes the proof of the 
 (\ref{e.distributional}) .

\begin{proof}
This is a consequence of the particular structure of a $+$tree $\Tree$, and the fact for $s\in \Tree$,  the distance of 
$\text{supp}(\widehat{\varphi_s})$ to $\omega_\Tree$ is approximately $\abs{\omega_s}$.   The Calder\'on Zygmund 
theory applies, and shows  that for any choice of signs $\varepsilon_s\in\{\pm1\}$, for $s\in \Tree$, 
 \begin{equation}  \label{e.cz} 
\ABs{ \Bigl\{ \sum_{s\in\Tree} \varepsilon_s \ip f,\varphi_s . \varphi_s >\lambda\Bigr\}}\lesssim\lambda^{-1}\norm f \abs{I_\Tree} \chi_{I_\Tree} .1. ,\qquad \lambda>0.
 \end{equation}
We apply this inequality for trees $\Tree \in \mathcal T_{\text{out}}$, and $f=\ind E  $.  
By taking the average over all choices of signs, we can conclude 
a distributional estimate on the square functions 
 \begin{equation} \label{e.zD}
\Delta_\Tree :=\Bigl[\sum_{s\in\Tree} \frac{ \abs{ \ip \ind E , \varphi_s .}^2 }{\abs {I_s}} \ind {I_s} \Bigr]^{1/2}  
 \end{equation}
Namely, that for each \tree+ $\Tree \subset \mathcal T_{\text{out}}$, 
 \begin{equation} \label{e.Gzl}
\abs{ \{ \Delta_\Tree >\lambda\} } {}\lesssim{}  \lambda^{-1} \abs{I_\Tree}\abs{F}^{-1}. 
 \end{equation}

As this inequality applies to all subtrees of $\Tree$, it can be strengthened.  (This is a reflection of the John Nirenberg inequality.) 
Fix the $+$tree $\Tree\subset\mathcal T_{\text{out}}$.  We wish to conclude that 
 \begin{equation}\label{e.sizeF} 
\abs{I_\Tree}^{-1}\int_{I_\Tree} \Delta_\Tree^2\; dx\lesssim\abs{F}^{-1}.
 \end{equation}
For a subset $\Tree'\subset\Tree$, let 
 \begin{equation*}
\operatorname{sh}(\Tree'):=\bigcup_{s\in\Tree'}I_s
 \end{equation*}
be the {\em  shadow of $\Tree'$.}  A shadow is not necessarily an interval. Define $\Delta_{\Tree'}$ as in (\ref{e.zD}) .  And finally set 
 \begin{equation}\label{e.GG}
G(\lambda)=\sup_{\Tree'\subset\Tree} \abs F \abs{\operatorname{sh}(\Tree')}^{-1} \abs{ \{ \Delta_{\Tree'} >\lambda\} }
 \end{equation}
Notice that (\ref{e.Gzl})  implies that $G(\lambda)\lesssim\lambda^{-1}$, for $\lambda>0$.  If we show that $G(\lambda)\lesssim{}\lambda^{-4}$, for $\lambda>1$, 
we can conclude (\ref{e.sizeF}) .   
In fact we can show that $G(\lambda)$ decays at an exponential squared rate, which is the optimal estimate.

Observe that (\ref{e.sizeF}) , implies that we have 
 \begin{equation*}
\frac{\abs{\ip \ind E ,\varphi_s . } } {\sqrt{ \abs I } }\le\lambda_0 <\infty. 
 \end{equation*}
Thus, the square functions we are considering $\Delta_\Tree $ can only take take incremental steps of a strictly bounded size.

For any $\lambda\ge\sqrt2\lambda_0$, let us bound $G(\sqrt 2\lambda)$.  
Fix $\Tree'$ achieving the supremum in the definition of $G(\sqrt2\lambda)$.  
Consider a somewhat smaller threshold,  namely  $\{ \Delta_{\Tree'}>\lambda\}$.   In order to proceed, consider a 
function $\tau\mid \operatorname{sh}  {\Tree}'\mapsto \mathbb R_+ $ such that 
 \begin{equation*}
\sum _{\substack{ s\in\Tree \\ \abs {I_s}\ge\tau(x) }}\frac{\abs{ \ip f,\varphi_s .}^2 } {{\abs {I_s}}} \ind {I_s}  (x)\ge\lambda^2. 
 \end{equation*}
In addition require that $\tau(x) $ is the smallest such function satisfying this condition.  It is the case that 
the sum above can  be no more than $\lambda^2+\lambda_0^2 $.  

Take $\Tree''\subset\Tree'$ to be the tree 
 \begin{equation*}
\Tree'':=\{s\in\Tree'\mid \abs{I_s }\le{}\tau(x),\ x\in I_s \}.
 \end{equation*}
 The point of these definitions is that 
  \begin{equation*}
 \Delta_{\Tree'}(x)\ge\sqrt2\lambda\quad\text{implies}\quad \Delta_{\Tree''}(x)\ge\sqrt{ \lambda^2-\lambda_0^2}.
  \end{equation*}
Therefore, 
 \begin{align*}
\abs{F}^{-1}{\operatorname{sh}(\Tree')} G(\sqrt2\lambda){}={}&
\abs{ \{\Delta_{\Tree'} >\sqrt2\lambda\} }
\\{}\le{}&  \abs{ \{\Delta_{\Tree''} >\sqrt{ \lambda^2-\lambda_0^2}\} }
\\{}\le{}&  \abs{F}^{-1}{\operatorname{sh}(\Tree'')} G(\sqrt{ \lambda^2-\lambda_0^2})
\\{}\le{}&  \abs{F}^{-1}{\operatorname{sh}(\Tree')} G(\lambda)G(\sqrt{ \lambda^2-\lambda_0^2}).
 \end{align*}
We conclude that $G(\sqrt2\lambda)\le{}G(\sqrt{ \lambda^2-\lambda_0^2})^2$.  

To conclude, we should in addition require that $\lambda_0 $ is so large that 
$
G(\kappa\lambda_0)\le{}\tfrac12 $, where  
 \begin{equation*}
\kappa:=\prod _{k=1}^\infty{} \sqrt{1-2^{-k} } . 
 \end{equation*}
An induction argument will then show that 
 \begin{equation}\label{e.exp}
G(2^{k/2}\lambda_0)\le{}G(\kappa\lambda_0)^{2^k}\le{}2^{-2^k},\qquad k\ge0,
 \end{equation}
which is the claimed exponential decay. Our proof of the Lemma is done.
\end{proof}

\subsection{Complements}  

 \begin{remark}\label{r.kp}   In the inequality  (\ref{e.distributional}) , one can show that the constants $K_p $ on the right hand side obey 
 \begin{equation*}
K_p\lesssim{} \frac{p^2}{p-1}
 \end{equation*}
 \end{remark}

 \begin{remark}\label{r.G} 
If it is the case that for some $0< \alpha<1 $, we have the inequality 
 \begin{equation*}
\sup_{\Tree'\subset\Tree} \abs{\operatorname {sh}(\Tree') }^{-1/\alpha} \norm \Delta_{\Tree'} .\alpha.<\infty 
 \end{equation*}
then, the stronger estimate below holds. 
 \begin{equation*}
\norm \Delta_\Tree .1. {}\lesssim{} \abs{\operatorname {sh}(\Tree) }. 
 \end{equation*}
 \end{remark}

 \begin{remark}\label{r.weakp}  In (\ref{e.distributional}) , we assert the restricted weak type inequality for $1<p<\infty$.  The weak type estimate 
for $2<p<\infty$ is in fact directly available.  That is, for $2<p<\infty$, and $f\in L^p$ of norm one,  
 \begin{equation*}
\abs{ \{ \mathcal C_N f >\lambda\}}\lesssim{} \lambda^{-p},\qquad \lambda>0.
 \end{equation*}
The key point is to take advantage of the fact that $f $ is locally square integrable.  A very brief sketch of the argument follows. 
 (1) It suffices to prove the inequality above for $\lambda=1$.  (2) Define 
$\Omega=\{ \operatorname M\abs{f}^2 >1\}$, and show that $\abs\Omega\lesssim1$.  (3) Define sums as in (\ref{e.inside})  and (\ref{e.outside}) , 
and control each term separately.   One will need to replace the Size Lemma as stated with \ref{r.ABove2}.
 \end{remark}

\section{Remarks}

 \begin{remark}\label{r.carleson} After L.~Carleson \cite{carleson} proved his theorem, R.~Hunt \cite{MR38:6296} extended the argument to 
$L^p$, for $1<p<\infty$.  A similar extension, in the Walsh Paley case, was done by Billiard \cite{MR36:599}, in the case 
of $L^2$, and Sj\"olin \cite{MR39:3222}, for all $1<p<\infty$.  The Carleson theorem has equivalent formulations on the groups 
$\mathbb R$, $\mathbb T$, and $\mathbb Z$.  The last case, of the integers, was explicitly discussed by M{\'a}t{\'e} \cite{MR39:701}. This paper was 
overlooked until recently.
 \end{remark} 

 \begin{remark}\label{r.feff}  C.~Fefferman \cite{fefferman} devised an alternate proof, which proved to be influential through it's 
 use of methods of analysis that used both time and frequency information in an operator theoretic fashion.  The proof 
 of Lacey and Thiele \cite{laceythielecarleson} presented here borrows several features of that proof. 
 The notion of tiles, and the partial order on tiles is due to C.~Fefferman \cite{fefferman}.  Likewise, the Density Lemma and the Tree Lemma, and the proof of the same, have clear antecedents in this paper. 
 \end{remark} 

 \begin{remark}\label{r.notseperated}  Those familiar with the Littlewood Paley theory know that it is very useful in decoupling the scales 
of operators, like those of the Hilbert transform.  The Carleson operator is, however, not one in which scales can be decoupled. 
This is another source of the interest in this Theorem. 
 \end{remark}

 \begin{remark}\label{r.equiv}  We present the proof of the Carleson theorem on the real line due to the presence of the dilation structure. 
 \end{remark}

 \begin{remark}\label{r.1}  We choose to express the Carleson operator in terms of the projection $\operatorname P_-$.  This operator is a linear combination of 
the identity operator and the Hilbert transform given by 
 \begin{equation*}
\operatorname H f(x):=\lim_{\epsilon\to0}\int_{\epsilon<\abs y<1/\epsilon} f(x-y)\;\frac {dy}y.
 \end{equation*}
Hence, an alternate form of the Carleson operator is 
 \begin{equation} \label{e.alternate}
\sup_N \abs{ \operatorname H \circ \modulate N f}=\sup_N\ABs{ \int  e^{iNy} f(x-y)\;\frac {dy}y}
 \end{equation}
This form is suggestive of other questions related to the Carleson Theorem, a point we rely upon below. 
 \end{remark}

 \begin{remark}\label{r.three} Despite the fact that Carleson's operator maps $L^2$ into itself, all three known proofs of Carleson's theorem establish 
the weak type bound on $L^2$.  The strong type bound must be deduced by interpolation.  On the other hand, the weak type bound is a known 
consequence of the pointwise convergence of Fourier series.  This was observed by A.~Calder\'on, as indicated by a footnote in \cite{zyg}, 
and is a corollary of a general observation of E.M.~Stein \cite{steinweak}.
 \end{remark}

 \begin{remark}\label{r.huntWtd}  Hunt and Young \cite{MR49:3419} have established a weighted estimate for the Carleson operator. Namely for a weight 
$w$ in the class $A^p$, the Carleson operator maps $L^p(w)$ into itself, for $1<p<\infty$.  The method of proof utilizes the known Carleson bound, and distribution inequalities for the Hilbert transform.
 \end{remark}

 \begin{remark}\label{r.P-}  The \Pcap.P-/ has a well known antecedent in a characterization of  (a constant times) the Hilbert transform 
as the unique operator $A$   such that $A$ is bounded operator on $L^2$ that commutes with dilation, is invariant under
dilations, $A^
2$ is the identity, but is not itself the identity.  See \cite{stein70s}.
 \end{remark}

 \begin{remark}\label{r.weaktypedual} 
The inequality (\ref{e.weaktypedual})  eschews all additional cancellations.  It shows that all the necessary cancellation properties are 
already encoded in the decomposition of the operator.   In addition the  combinatorial model of the 
Carleson operator is in fact unconditionally convergent in $s\in\mathcal T$.  This turns out to be extremely useful 
fact in the course of the proof: one is free 
to 
group the tiles in anyway that one likes.
 \end{remark}

 \begin{remark}\label{r.rubio<}  The Size Lemma should be compared to Rubio de Francia's extension of the classical Littlewood Paley inequality 
\cite{rubio}.  Also see the author's recent survey of this theorem \cite{laceyrubio}. 
 \end{remark}

 \begin{remark}\label{r.+}  A $+$tree $\Tree$ is a familiar object.  Aside from a modulation by $c(\omega_\Tree)$, it shares most
of the properties associated with sums of wavelets.  In particular, if $0\in\omega_\Tree$, note that 
 \begin{equation*}
\size \Tree {}\simeq{} \NOrm \sum_{s\in\Tree} \langle f,\varphi_s \rangle  \varphi_s .BMO.
 \end{equation*}
where the last norm is the $BMO$ norm.
 \end{remark}

 \begin{remark}\label{r.tree}  The key instance of the Tree Lemma is that of a $+$tree.  This case corresponds to a particular maximal function 
applied to a function associated to the tree. 
It is this point at which the supremum of Carleson's theorem is controlled by a much tamer supremum:  The one in the ordinary maximal function.
 \end{remark}

 \begin{remark}\label{r.size} The statement and proof of the size lemma \l.size/ replaces the initial arguments of this type that are in 
Lacey and Thiele \cite{MR99b:42014}.  This argument has proven to be very flexible in it's application. And, in some instances 
it produces sharp estimates, as explained by Barrioneuvo and Lacey \cite{MR1955263}.
 \end{remark}

 \begin{remark}\label{r.gabor}  The set of functions $\mathcal T_k:=\{s\in\mathcal T\mid \abs{I_s}=2^k\}$ is an example of a Gabor basis.  For 
appropriate choice of $\varphi$, the operator  
 \begin{equation*}
\operatorname A_kf:=\sum_{\substack{s\in\mathcal T \\ \abs{I_s}=2^k}}  \ip f,\varphi_s. \varphi_s 
 \end{equation*}
is in fact the identity operator.  See the survey of I.~Daubechies \cite{MR93e:42045}. 
 \end{remark}

 \section{Complements and Extensions}

\subsection{Equivalent Formulations of Carleson's Theorem}

The Fourier transform has a formulation on each of the Euclidean groups 
$\mathbb R$, $\mathbb Z$ and $\mathbb T$. 
Carleson's original proof worked on $\mathbb T$. 
Fefferman's proof translates very easily to $\mathbb R$.  M{\'a}t{\'e}
\cite{MR39:701} extended Carleson's proof to $\mathbb Z$.
   Each of the statements of the theorem can be stated in terms of a maximal Fourier multiplier theorem, and we 
   have stated it as such in this paper. 
Inequalities for such operators can be transferred  between these three Euclidean groups, and was done so by 
P.~Auscher and M.J.~Carro \cite{MR94b:42007}.  

Transference has also been studied in the bilinear setting.  See articles by 
Fan and Sato, as well as Blasco and Villarroya 
\cites{MR2037006,MR1808390}.

\subsection{Fourier Series Near $L^1$, Part 1}  \label{s.near1}

The point of issue here is the  determination of that integrability class which 
guarantees the pointwise convergence of Fourier series. The natural setting for these questions is the 
unit circle $\mathbb T=[0,1]$, and the partial Fourier sums 
 \begin{equation*}
S_N f(x)=\sum_{\abs n\le N} \widehat f (n) e^{2\pi inx},\qquad \widehat f (n)=\int_0^1 f(x)e^{-inx}\; dx
 \end{equation*}
  In the positive direction, one seeks the 
``smallest'' function $\psi$ such that if $\int_\mathbb T \psi(f)\; dx<\infty$, then the Fourier series of $f$ converge pointwise.

 N.~Antonov \cite{MR97h:42005} has found the best result to date, 
 
\begin{theorem} \label{t.antonov}
For all functions $f\in L(\log L)(\log\log\log L)(\mathbb T)$, the partial Fourier series of $f$ converges pointwise to $f$. 
\end{theorem}

This extends the result of P.~Sj\"olin \cites{MR39:3222,MR49:998}, who had the result above, but with a double log where there is a triple log above. 
  Arias de Reyna \cites{MR1906800,MR2002k:42009} has noted an extension of this Theorem, in that one can define a rearrangement invariant Banach space $B$, so that pointwise convergence holds for all $f\in B$, and $B$ contains $L(\log L)(\log\log\log L)$.

The method of proof  takes as it's starting point, the distributional estimate of (\ref{e.best}) .
  One seeks to ``extrapolate'' these inequalities to the setting of the Theorem above   and Antonov nicely exploits the 
  explicit nature of the kernels involved in this maximal operator. 
  Also see the work of P.~Sj\"olin and  F.~Soria \cite{MR2014553}  who demonstrate that Antonov's approach 
  extends to other maximal operator questions.

\subsection{Fourier Series Near $L^1$, Part 2}

In the negative direction, 
A.N.~Kolmogorov's fundamental example  \cites{kolo1,kolo2} of an integrable function with pointwise
divergent Fourier series admits a strengthening to the following statement, as obtained by K\"orner \cite{MR83h:42010}.

\begin{theorem} \label{t.kolo}  For all  $\psi(x)=o(\log \log x)$, there   is a function $f\mid [0,2\pi]\to\mathbb R$ with divergent Fourier series, 
and $\int \abs f \psi(f)\; dx<\infty$. 
\end{theorem}

The underlying method of proof was, in some essential way, unsurpassed until quite recently, when S.~Konyagin \cites{MR2000h:42002,konyagin}  proved 

\begin{theorem} \label{t.konyagin} The previous Theorem holds assuming only 
 \begin{equation} \label{e.kon} 
\psi(x)=o\Biggl( \sqrt{ \frac{ {\log x}}{ \log\log x}}\Biggr)
 \end{equation}
\end{theorem}

There is a related question, on the growth of partial sums of Fourier series of integrable functions.  G.H.~Hardy \cite{hardy} showed that for 
integrable functions $f$, one has $S_nf=o(\log n)$ a.e., and asked if this is the best possible estimate.  This question is still open, with the best result from below following from Konyagin's example.  With $\psi$ as in (\ref{e.kon}) , there is an $f\in
 L^1(\mathbb T
)$ with 
 \begin{equation*}
\limsup_{n\to\infty} \frac{S_n f}{\psi(n)}=\infty \qquad \text{for all $x\in\mathbb T$.}
 \end{equation*}
 Bochkarev \cite{MR1993585} has  a very slight strengthening of this result for Walsh series.

\bigskip 

Let $\delta_t$ denote the Dirac point mass at $t\in\mathbb T$.  The method of proof is to construct measures 
 \begin{equation*}
\mu=K^{-1}\sum_{k=1}^K \delta_{t_k}
 \end{equation*}
a set $E\subset \mathbb T$ with measure at least $1/4$, and choices of integers $N$ for which 
 \begin{equation*}
\sup_{n<N} \abs{ S_n \mu(x)}\ge \psi(N),\qquad x\in E.
 \end{equation*}
Kolmogorov's example  consists of uniformly distributed point masses, whereas Konyagin's example consists of point masses 
that have a distribution reminiscent of a Cantor set.  

\subsection{Probabilistic Series} 

It is of interest from the point of view of probability and ergodic theory, to consider the version of the Hilbert transform 
and Carleson theorem that arises from the integers.  Here, we consider the probabilistic versions.  Let $ X_k $ be independent 
and identically distributed copies of a mean zero random variable $X $.   The question is  if the sum 
 \begin{equation*}
\sum_{k=1}^ \infty \frac{X_k} k 
 \end{equation*}
converges a.s.  Without additional assumption on the distribution of  $X $, a necessary and sufficient condition is that 
$\mathbb E X \log (2+\abs X) <{} \infty $.  One direction of this is in \cite{mz}.   If however  $X $ is assumed to be symmetric, 
integrability is necessary and sufficient.   This addresses the issue of the Hilbert transform.  

Carleson's theorem, in this language, concerns the convergence of the series 
 \begin{equation*}
Y(t):=\sum_{k=1}^ \infty \frac{X_k} k  e^{2\pi i k t}   \qquad\text{for all $t\in \mathbb T $.} 
 \end{equation*}
The role of the quantifiers should be emphasized.  Convergence holds for all $t\in \mathbb T $,  on a set of full probability. 
 Given this, abstract results on $0$\ndash$1 $ Laws assure us that if the series converges for all $t $,  off of a single set 
of probability zero, then the limiting function is continuous with probability one.  
The paper of  Talagrand \cite{MR96f:60068} gives necessary and sufficient conditions for the convergence of this series.

 \begin{theorem}\label{t.talagrand}  Let $X_k$ be independent identically distributed copies of a mean zero symmetric random variable $X$.  
$X\in L(\log \log L)$ iff the series  $Y(t) $ converges to a continuous function on $\mathbb T$ almost surely. 
 \end{theorem}

The assumption of symmetry should be added to the statement of the Theorem in 
\cite{MR96f:60068}.   Cuziak and Li \cite{MR82f:60083} provide an example of a non symmetric mean zero $X \in L (\log\log L) $ for 
which the series $Y(t) $ is divergent.    This series is a borderline series in that it just falls out of the scope of the powerful 
theory of Marcus and Pisier \cite{MR86b:60069} on random Fourier series.  

My thanks to several people who provided me with some references for this section.  They are James Campbell, Ciprian Demeter, 
Michael Lin, and Anthony Quas.

\subsection{The Wiener--Wintner Question}  \label{s.ww}

A formulation of Carleson's operator on $\mathbb Z$ is 
$$
C_\mathbb Z f(j):=\sup_\tau\sup_N\ABs{ \sum_{0<\abs k<N} f(j-k)\frac{e^{i\tau k}}k}.
$$
See M{\'a}t{\'e} \cite{MR39:701}.  Unaware of this work which followed soon after 
Carleson,  J.~Campbell and Petersen \cite{MR89i:28008} considered this operator on $\ell^2$, 
with equivalence in $\ell^p$ established by  Assani and Petersen \cites{MR89i:28008,MR92g:28025}.  Also see Assani, Petersen and White \cite{MR92m:28021}, for these and other equivalences.  The latter  authors  had additional motivations from dynamical sy

stems, which we turn to now.

A.~Calder\'on \cite{MR37:2939} observed that 
inequalities for  operators on $\mathbb Z$ which commute with translation  can be transferred to
discrete dynamical systems.  Let $(X,\mu)$ be a probability space, and $T\mid X\to X$ a map which preserves 
$\mu$ measure.  Thus, $\mu(T^{-1}A)=\mu(A)$ for all measurable $A\subset X$.  A Carleson operator on $(X,\mu,T)$ is 
$$
\operatorname C_{\text{mps}} f(x):=\sup_\tau\sup_N\ABs{ \sum_{0<\abs k<N} f(T^kx)\frac{e^{i\tau k}}k}
$$
And it is a consequence of Calder\'on's observation and Carleson's theorem that this operator is bounded on $L^2(X)$.  

There is however a curious point that distinguishes this case from the other settings of Euclidean groups.  It is the case that one has pointwise convergence of 
$$
\lim_{N\to\infty} \sum_{0<\abs k<N} f(T^kx)\frac{e^{i\tau k}}k \quad\text{exists for all $\tau$}
$$
holding for almost every $x\in X$?  The boundedness of the maximal function $C_{\text{mps}}$ shows that this would hold on a closed set in $L^2(X)$.  The missing ingredient is the dense class for which the convergence above holds.
  Unlike the setting of Euclidean groups, there is no natural dense class.

This conjecture was posed by Campbell and Petersen \cite{MR89i:28008}.

 \begin{conjecture}\label{j.ww}  
For all measure preserving systems $(X,\mu,T)$, and all $f\in L^2(X)$, we have the following:
$$
\mu\Biggl\{ \lim_{N\to\infty} \sum_{0<\abs k<N} f(T^kx)\frac{e^{i\tau k}}k \quad\text{exists for all $\tau$} \Biggr\}=1.
$$
 \end{conjecture}

\medskip 

A theorem of N.~Wiener and A.~Wintner \cite{MR4:15a} provides a classical motivation of this question.  This theorem concerns the same phenomena, but with the discrete Hilbert transform replaced by the averages.

\begin{theorem}  \label{t.ww}  
For all measure preserving systems $(X,\mu,T)$, and all $f\in L^2(X)$, we have the following:
$$
\mu\Biggl\{ \lim_{N\to\infty} N^{-1} \sum_{k=0}^{N-1} f(T^kx){e^{i\tau k}} \quad\text{exists for all $\tau$} \Biggr\}=1.
$$
\end{theorem}

This theorem admits a simple proof.  And note that this theorem trivially supplies a dense class in all $L^p$ spaces, $1\le{}p<\infty$.  

The Wiener--Wintner theorem has several interesting variants, 
for which one can phrase related questions by replacing averages by Hilbert transforms. 
As far as is known to us, none of these questions is answered. 

\medskip 

An attractive theorem proved by E.~Lesigne \cites{MR91h:28016,MR95e:28012} is 

\begin{theorem} For any measure preserving system $(X,\mu,T)$ and all integrable functions $f$, 
there is a subset $X_f\subset X$ of full measure so that for all $x\in X_f$, all polynomials $p$, and all 
1--periodic functions $\phi$, the limit below exists:

$$
\lim_{N\to\infty}N^{-1}\sum_{n=1}^N \phi(p(n))f(T^nx)
$$
\end{theorem}

Extending this theorem to the Hilbert transform would be an extraordinary accomplishment, whereas if one replaced the discrete dynamical system by flows, it could be that the corresponding result for the Hilbert transform might be within reach.   

In connection to this, Arkihpov and Oskolkov \cite{MR89a:42010} have proved the following theorem.

 \begin{theorem}\label{t.arkihpov}  For all integers $d$, 
 \begin{equation*}
\sup_{\operatorname{deg}(p)=d} \ABs{ \sum_{n\not=0}\frac{ e^{i p(n)} }n }<\infty.
 \end{equation*}
with the supremum formed over all polynomials of degree $d$.
 \end{theorem} 

This is a far more subtle fact than the continuous analog stated in (\ref{e.swp}) .  Arkihpov and Oskolkov use the Hardy Littlewood 
Circle method of exponential sums, with the refinements of Vinogradov.  See \cites{MR89a:42010,MR99j:42004,MR94m:42016}.  

By  Plancherel,  this theorem shows that for a polynomial $p$ which maps the integers to the integers,
 the operators on the integers given by 
 \begin{equation*}
T_p f(j)=\sum_{n\not=0}f(x-p(n))
 \end{equation*}
is a bounded operator on $\ell^2(\mathbb Z)$.

E.M.~Stein and S.~Wainger have established $\ell^2$ mapping properties for certain 
Radon transforms \cites{MR2000j:42027,MR2001g:42038}.

\subsection{E.M.~Stein's Maximal Function}   \label{s.stein} 

A prominent theme of the research of E.M.~Stein and S.~Wainger concerns oscillatory integrals, with polynomial phases.  It turns out to be of interest to determine what characteristics of the 
polynomial govern allied analytic quantities.  In many instances, this characteristic is just 
the the degree of the polynomial.  For instance, the following is a corollary to a 
Theorem of Stein and Wainger from 
1970 \cite{MR42:904}.  Namely, that one has a bound 
 \begin{equation} \label{e.swp}
\sup_{\operatorname{deg}(P)=d} \ABs{ \int e^{ip(y)}\;\frac{dy}y}\lesssim{}1, \qquad d=1,2,\ldots
 \end{equation}

A conjecture of Stein's concerns an extension of Carleson's maximal operator to one in which one forms a supremum over 
all polynomial choices of phase with a fixed degree.  Thus, 

 \begin{conjecture}\label{j.stein}  For each integer $d$, the maximal function below maps $L^p$ into itself for $1<p<\infty$.  
 \begin{equation*}
\mathcal C_df(x)=\sup_{\operatorname{deg}(P)=d} \ABs{ \int e^{ip(y)} f(x-y)\;\frac{dy}y}
 \end{equation*}
 \end{conjecture}

Note that the case of $d=1$ corresponds to Carleson's theorem.   Let us set $\mathcal C'_d$ to be the maximal operator above, but with the 
the restriction that the polynomials $p$ do not have a linear term.  It is useful to make this distinction, as it is the linear terms 
that are intertwined with the Fourier transform.

E.M.~Stein   \cite{MR96i:42013} considered the purely quadratic terms, and showed that $\mathcal C'_2$ maps $L^p$ into itself
for all $1<p<\infty$.
The essence of the matter is the bound on $L^2$, and there his argument is a variant on  the method of $\operatorname T
\operatorname T^*$, emphasizing 
a frequency decomposition of the operator.    Stein and Wainger \cite{MR2002k:42038} have  proved that $\mathcal C'_d$
is bounded on all $L^p$'s, for all $d\ge2$. Again 
the $L^2$ case is decisive and the argument is  an application of the $\operatorname T \operatorname T^*$ method,
but with a spatial decomposition of the operator.

Let us comment in a little more detail about how these results are proved.  If, for the moment, one consider a fixed polynomial 
$P(y) $, and the oscillatory integral 
 \begin{equation} \label{e.TP} 
T_P f(x):=\int e^{iP(y) } f(x-y) \frac{dy}y .
 \end{equation}
One may utilize the scale invariance of the the Hilbert transform kernel to change variables.  With the correct change of variables, 
one may assume that the polynomial $P(y)=\sum _{j=1}^d a_j y^j $ satisfies $\sum \abs{a_j}=1 $.   Then, it is evident that 
for $\abs y <1 $, say, that the integral above is well approximated by a truncation of the Hilbert transform.   Thus, it is those 
scales of the operator larger than $1 $ that must be controlled.  

It is a consequence of the van der Corput estimates that some additional decay can be obtained from these terms.  
In particular, one has this estimate.   To set 
notation,  in the one dimensional case only, set  
 \begin{equation*} 
P_{\vec a}(x)=a_d x^d+\cdots+{} a_1 x,\qquad  \vec a=(a_d,\ldots,a_1).
 \end{equation*}

 \begin{lemma}\label{l.vdc}  Let $\chi $ be a smooth bump function. Then we have the estimate 
 \begin{equation*}
\NOrm \widehat{\operatorname e ^{iP_{\vec a} (y) } \chi(y) }(\xi). \infty. {}\lesssim{}
(1+\norm \vec a .1.)^{-1/d}. 
 \end{equation*}  
In particular, by the Plancherel identity, we have the estimate 
 \begin{equation} \label{e.vdc} 
\Norm [\operatorname e ^{iP_{\vec a} (y) } \chi(y)]*f(x) .2. {}\lesssim{} (1+\norm \vec a .1.)^{-1/d} \norm f.2. .
 \end{equation}
 \end{lemma}

Notice that these estimates are better than the trivial ones.  And that the second estimate can be interpolated to 
obtain a range of $L^p $ inequalities for which one has decay, with a rate that depends upon  the degree and the $L^p $ space in 
question. 

In a discussion of the extensions of this principle in  for example \cites{MR2002k:42038,MR96i:42013}, one establishes 
appropriate extensions of this last lemma, always seeking some additional decay that arises from the polynomial.  
For instance, in \cite{MR2002k:42038}, Stein and Wainger prove a far reaching extension of this principle.

 \begin{lemma}\label{l.sw}   There is a constant $\delta>0 $, depending only on the degree $d $, so that we have the estimate 
$\norm S_ \lambda .2. {}\lesssim{}(1+\lambda)^{-\delta} $, for all $\lambda>0 $, where 
 \begin{equation*}
\operatorname  S _{\lambda}f(x):={} 
\sup_{\substack{\norm\vec a.1.\ge\lambda\\ a_1=0 }}  \sup _{t>0 } 
\abs{ [\operatorname  {Dil}_t^1 \operatorname e ^{i P_{\vec a}(\cdot )} \chi ]*f }
 \end{equation*}
It is essential in this supremum be formed over polynomials $P_{\vec a} $ which do not have a linear term.   
 \end{lemma} 

A.~Ionescu has pointed out that this Lemma is not true with the linear 
term included, even in the case  of second degree polynomials.
The example, which we will see again below, 
begins by taking a function $f (x) $, and replacing it 
by the function $g(x)=\operatorname e ^{i \lambda x^2} f(x)$. 
Then, in the  supremum defining $S_ \lambda $ above, 
take the dilation parameter to be $t=1 $, 
and the polynomial to be $P(y)=-y^2+2xy $.  Note that as we
are taking a supremum, we can in particular take a polynomial 
that depends upon $x $.  

In this example, the modulation of $f $ by ``chirp''{}  is then canceled out by the choice of $P $.  There is no decay in the 
estimate. This estimate is special to the case of the second power, so it is natural to guess that it 
plays a distinguished role in these considerations.

This also  points out an error in the author's paper \cite{MR1949031}.  (The error enters in specifically at the equation 
(2.9). The phase plane analysis of that paper might yet find some use.)
At this point, the resolution of Stein's conjecture is not settled.  And it appears that a positive bound of the 
operator $\mathcal C_2 $ will in particular require a novel phase plane analysis with quadratic phase. 
This  should be compared to the notion of \emph{degeneracy} in \s.multi/ below.

\subsection{Fourier Series in Two Dimensions}

In this section we extend the Fourier transform to functions of the plane 
 \begin{equation*}
\widehat f(\xi)=\int f(x)e^{ix\cdot \xi}\; dx
 \end{equation*}
where $\xi=(\xi_1,\xi_2)$, and $x=(x_1,x_2)$.  The possible extensions of Carleson's theorem 
to the two dimensional setting are numerous. The state of our knowledge is not so great. 

\subsubsection{Fourier Series in Two Dimensions, Part I} 

Consider the pointwise convergence of the Fourier sums in the plane given by 
 \begin{equation*}
\int_{tP} \widehat f(\xi)e^{ix\cdot\xi}\; dx 
 \end{equation*}
Here, $P$ is   polygon with finitely many sides in the plane, with the origin in it's interior.  
  Proving the pointwise convergence of 
these averages is controlled by the maximal function 
 \begin{equation*}
\mathcal C_P f {}:={} \sup_{t>0} \ABs{ \int_{tP} \widehat f(\xi)e^{ix\cdot\xi}\; dx }
 \end{equation*}
C.~Fefferman \cite{MR55:8682} has observed that this maximal operator can be controlled by a sum   of 
operators which are equivalent to the Carleson operator. 

For simplicity, we just assume that the polygon is the unit square.  And let 
 \begin{equation*}
Qf=\int_0^\infty \int_{-\xi_1}^{\xi_1} \widehat f(\xi_1,\xi_2) e^{i (x_1\xi_1+x_2\xi_2)} \;d\xi_2d\xi_1
 \end{equation*}
This is the Fourier projection of $f$ onto the sector swept out by the right hand side of the square. 
Notice that $\mathcal C_P\circ Q$ is the one dimensional Carleson operator applied in the first coordinate.  

Thus, $\mathcal C_P$ is dominated by a sum of terms which are obtained from  the one dimensional Carleson operator, 
and so $\mathcal C_P$ maps $L^p$ into itself for $1<p<\infty$.  This argument works for any polygon with a finite number of 
sides.  While we have stressed the two dimensional aspect of this argument, it also works in any dimension.

\medskip  

Nevertheless, it is of some interest to consider maximal operators of the form 
 \begin{equation*}
\sup_{\xi\in\mathbb R^2}\ABs{ \int f(x-y) e^{i\xi \cdot y} K(y)\; dy }                                                                       
 \end{equation*}
where $K$ is a Calder\'on Zygmund kernel.  This is the question addressed by P.~Sj\"olin \cite{MR55:8682}, and more recently by 
Sj\"olin and Prestini \cite{MR2001m:42017} and Grafakos, Tao and Terwilleger \cite{MR2031458}. 

\subsection{Fourier Series in Two Dimensions, Part II} 
What other methods can be used to sum Fourier series in the plane?  One method that comes to mind is over arbitrary rectangles.  That is, 
one considers the maximal operator
 \begin{equation*}
\mathcal R f(x):=\sup_\omega \ABs{ \int_\omega \widehat f(\xi) e^{i\xi\cdot x}\; d\xi } 
 \end{equation*}
The supremum is formed over arbitrary rectangles $\omega$ with center at the origin.
C.~Fefferman \cite{MR43:5251} has shown however that this is a  badly behaved operator. 

 \begin{proposition}\label{p.multiple}  There is a bounded, compactly supported function $f$ for which $\mathcal R f=\infty$ on 
a set of positive measure. 
 \end{proposition} 

This maximal operator has an alternate formulation, see (\ref{e.alternate}) , as 
 \begin{equation} \label{e.recCar}
\sup_{\alpha,\beta}\ABs{ \int\!\!\int f(x-x',y-y') e^{i(\alpha x'+\beta y')} \; \frac  {dx'}{x'}\frac  {dy'}{y'} }
 \end{equation}
The example of C.~Fefferman is  a sum of terms  $f_\lambda(x,y)=e^{i\lambda x y}\chi (x,y)$, where $\lambda>3$, and  $\chi$ is a 
smooth  bump function satisfying e.g.~$\ind {[-1,1]^2}  \le\chi\le\ind {[-2,2]}  $.  The key observation in the construction 
of the example is 

 \begin{lemma}\label{l.feff2}  We have the estimate 
 \begin{equation*}
\mathcal R f_\lambda (x,y) \gtrsim\log \lambda 
 \end{equation*}
for $ (x,y)\in [-\frac12,\frac12]^2$.
 \end{lemma} 

\begin{proof}
In the supremum over $\alpha$ and $\beta$ in the definition of $\mathcal R$, let $\alpha=\lambda y$ and $\beta=\lambda x$, and consider 
 \begin{align*}
R(x,y)={}&\int \Bigl[ \int f_\lambda(x-x',y-y')e^{i\lambda (x'y+xy')} \; \frac  {dx'}{x'}\Bigr]\frac  {dy'}{y'} 
\cr
{}={}&e^{i\lambda xy}\int \Bigl[ \int e^{i\lambda x'y'} \chi(x-x',y-y')\; \frac  {dx'}{x'}\Bigr]\frac  {dy'}{y'} 
 \end{align*}
The inside integral in the brackets admits these two estimates for all $(x,y)\in[-\frac12,\frac12]^2$.  
 \begin{align*} 
I(x,y,y')={}&\Bigl[ \int e^{i\lambda x'y'} \chi(x-x',y-y')\; \frac  {dx'}{x'}\Bigr]
\cr
{}={}&\begin{cases} 
c \operatorname{sign}(\lambda y')+O( \lambda\abs{y'})^{-1} & \cr
O(1+ \lambda \abs{y'}) & \cr
\end{cases}
 \end{align*}
$c$ is a non--zero constant. 
Both of these estimates are well--known. 

Then, we should estimate 
 \begin{equation*}
R(x,y)=\int_{\abs {y-y'}<1/\lambda}   I(x,y,y')\frac  {dy'}{y'} {}+\int_{\abs {y-y'}>1/\lambda}   I(x,y,y')\frac  {dy'}{y'} .
 \end{equation*}
The first term on the right is no more than $O(1)$, and the second term is $\gtrsim\log 1/\lambda$.
\end{proof}

This example shows that there are bounded functions $f$ for which 
 \begin{equation*}
\sup_{\abs\alpha,\abs\beta<N}\ABs{ \int\!\!\int f(x-x',y-y') e^{i(\alpha x'+\beta y')} \; \frac  {dx'}{x'}\frac  {dy'}{y'} }\gtrsim\log N. 
 \end{equation*}
It might be of interest to know if this estimate is best possible.

The integrals  in (\ref{e.recCar})  are  singular integrals in the product setting.  
There is as of yet no positive results relating to Carleson's theorem in a product setting.  

\subsubsection{Fourier Series in Two Dimensions, Part III}  \label{s.spherical}

The exponential $e^{i\xi\cdot x}$ is an eigenfunction of $-\Delta$, the positive Laplacian, with eigenvalue $\abs{\xi}^2$.  It would 
be appropriate to sum Fourier series according to this quantity.  This concerns the operator of Fourier restriction to the unit 
disc 
 \begin{equation*}
Tf:=\int_{\abs \xi<1}\widehat f (\xi) e^{i\xi \cdot x}\; d\xi 
 \end{equation*}
It is a famous result of C.~Fefferman \cite{MR45:5661} that 
$T$ is  bounded on $L^2(\mathbb R^2)$ iff $p=2$.  The fundamental reason for this is the existence of the Besicovitch set, a set contained in a large ball, that contains a unit line segment in each dir
ection, but has Lebesgue measure one.    The relevance of this set is indicated by the observation that the restriction of $T$ to  very small disc placed on the boundary of the disc is well approximated by 
a projection onto a half space.  Such a projection is a one dimensional Fourier projection performed in the normal direction to the disc.  And the normal directions can point in arbitrary directions. 
An extension of Fefferman's argument shows that the Fourier restriction to any smooth set with a curved boundary can only be a bounded operator on $L^2$.  

Nevertheless, the question of summability in the plane remains open.  Namely, 

\begin{conjecture} \label{q.spherical}  Is it the case that the maximal operator below maps $L^2(\mathbb R^2)$ into weak $L^2(\mathbb R^2)$? 
 \begin{equation*}
\sup_{r>0} \ABs{ \int_{\abs \xi<r}\widehat f (\xi) e^{i\xi \cdot x}\; d\xi }
 \end{equation*}
\end{conjecture}

In the known proofs of Carleson's theorem, the truncations of singular integrals plays a distinguished role.
In the proof we have presented, this is seen in the Tree Lemma, also cf.~\r.tree/.  In view of this, it interesting to suppose that if this
conjecture is true, what could play a role similar to the Tree Lemma.  It appears to be this. 

\begin{conjecture}\label{q.sphericallacunary}
   Is it the case that the maximal operator below maps $L^2(\mathbb R^2)$ into weak $L^2(\mathbb R^2)$? 
 \begin{equation*}
\sup_{k\in\mathbb N} \ABs{ \int_{\abs \xi<1+2^{-k}}\widehat f (\xi) e^{i\xi \cdot x}\; d\xi }
 \end{equation*}
 
 \end{conjecture} 

It is conceivable that a positive answer here could lead to a proof of spherical convergence of Fourier series. 
A recent relevant paper is \cite{MR2003k:42027} by Carbery, Georges, Marletta and Thiele.

\subsubsection{Fourier Series in Two Dimensions, Part IV}  

In order to bridge the gaps between Parts I and III, the following question comes to mind.  Is there 
a polygon with infinitely many sides which one could sum Fourier series with respect to?

G.~Mockenhoupt pointed out to the author that there is a natural first choice for $P$.  It is a polygon $  P_{\text{lac}}$ which in the first quadrant has vertices at the points $e^{\pi i2^{-k}}$ for $k\in\mathbb N$. Call this the \emph{lacunary sided pol
ygon}.

 It is a  fact due to Cordoba and R.~Fefferman \cite{MR55:6096} that the lacunary sided polygon is a bounded $L^p$ multiplier, for all $1<p<\infty$. 
 That is the operator below maps $L^p$ into itself for $1<p<\infty$. 
  \begin{equation*}
 \operatorname T_{\text{lac}} f(x)=\int_{P_{\text{\rm lac}}} \widehat f(\xi)e^{i\xi\cdot x}\; d\xi 
  \end{equation*}
  This fact is in turn linked to the boundedness of the maximal function in a lacunary set of directions:
$$
\operatorname M_{\text{lac}}f(x)=\sup_{k\in\mathbb N}\sup_{t>0}
(2t)^{-1}\int_{-t}^t \abs{f(x-u(1,2^{-k}))}\; du
$$
Note that this is a one dimensional maximal function computed in a set of directions in the plane that, in a strong sense, is zero dimensional. 
Let us state as a conjecture.
 
 \begin{conjecture}\label{j.lacsided}  For $2\le{}p<\infty$, the maximal function below  maps $L^p(\mathbb R^2)$ into itself.
 \begin{equation*}
\sup_{t>0} \ABs{ \int_{t P_{\text{\rm lac}}} \widehat f(\xi)e^{i\xi\cdot x}\; d\xi }
  \end{equation*}
 \end{conjecture} 

Even a restricted  version of this conjecture remains quite challenging.  In analogy to \q.sphericallacunary/

 \begin{conjecture}\label{j.lacsidedlac}  For $2\le{}p<\infty$, the maximal function below  maps $L^p(\mathbb R^2)$ into itself.
 \begin{equation*}
\sup_{k\in\mathbb N} \ABs{ \int_{(1+2^{-k}) P_{\text{\rm lac}}} \widehat f(\xi)e^{i\xi\cdot x}\; d\xi }
  \end{equation*}
 \end{conjecture}

Another question, with a somewhat more quantitative focus,  considers uniform polygons with $N$ sides, but then seek norm bounds on these two maximal operators, on $L^2$ say,    
which grow logarithmically in $N$.  We do not have a good conjecture as to
the correct order of growth of these constants.   If one could prove that the bounds where independent of $N$, then
the spherical summation conjecture would be a consequence.  

\subsubsection{Fourier Series in Two Dimensions, Part V}  

Again in the plane, consider the Fourier restriction to the semi infinite rectangles 
 \begin{equation*}
\operatorname S _{n} f(x):=\int_{-\infty}^n\int_{-\infty}^{n^2} \widehat f (\xi) e^{i \xi\cdot x }\; d \xi,\qquad n>0.
 \end{equation*}
Thus, we are restricting to a semi infinite rectangle with vertex on a parabola.  If we consider the maximal operator 
 \begin{equation*}
\operatorname S f:=\sup_{n>0} \abs{\operatorname S_n f }.
 \end{equation*}
The question then concerns the $L^p $ boundedness of this operator on $L^p $ spaces.  

To this end, we remark that a very nice argument of C.~Fefferman \cite{MR55:8682} shows that the

\subsection{The Bilinear Hilbert Transforms} 

The bilinear Hilbert transforms are given by 
 \begin{equation*}
\operatorname H_\alpha (f,g)(x):=\int f(x-\alpha y)g(x-y)\;\frac{dy}y ,\qquad \alpha\in\mathbb R,
 \end{equation*}
with the convention that $\operatorname H_0(f,g)=f\operatorname  Hg$, and $\operatorname H_\infty(f,g)=(\operatorname Hf)g$.  A third degenerate value is $\alpha=1$.  

These transforms commute with appropriate joint translations of $f$ and $g$, and dilations of $f$ and $g$.  
They are related to Carleson's theorem through the observation that for $\alpha\not\in\{0,1,\infty\}$, 
$\operatorname H_\alpha$ enjoys an invariance property with respect to modulation.  Namely, 
 \begin{equation*}
\operatorname H_\alpha( \modulate \beta f, \modulate {-\alpha\beta} g)=\modulate {\beta-\alpha\beta} \operatorname H_\alpha(f,g).
 \end{equation*}
That is, the bilinear Hilbert transforms share the essential characteristics of the Carleson operator.  

It was the study of these transforms that lead Lacey and Thiele to the proof of Carleson's theorem presented here. 
The bilinear Hilbert transforms are themselves interesting objects, with surprising properties.  Indeed, it is natural 
to ask what $L^p$ mapping properties are enjoyed by these transforms.  Note that in the integral, the term $dy/y$ is dimensionless, 
so that the $L^p$ mapping properties should be those of H\"older's inequality.  Thus, $\operatorname H_\alpha$ should map $L^2\times L^2$  into $L^1$.  
Note that this is false in the degenerate case of $\alpha=1$, as the Hilbert transform does not preserve $L^1$.  
Nevertheless, this   was conjectured by A.~Calder\'on in the non--degenerate cases. 

 See \cites{MR2000d:42003,MR99e:42013,MR99b:42014,MR98e:44001} for a proof of this theorem. 

 \begin{theorem}\label{t.bht}  For $1<p,q\le\infty$, if  $0< 1/r=1/p+1/q<3/2$, and $\alpha\not\in\{0,1,\infty\}$, then 
 \begin{equation*}
\norm \operatorname H_\alpha (f,g).r.\lesssim{}\norm f.p.\norm g.q.
 \end{equation*}
 \end{theorem}

We should mention that in a certain sense the proof of this theorem is  easier than that for Carleson's operator.  The proof outlined in 
\cite{laceyicm} contains the notions of tiles, trees, and size.  But the estimate that corresponds to the tree lemma is a triviality.  
The reason for this gain in simplicity is that there is no need for a mechanism to control a supremum.

The subject of multilinear operators with modulation invariance has  inspired a large number of results, and is worthy of survey on 
it alone.  We refer the reader to  C.~Thiele's article   \cite{thieleicm} for a survey of recent activity in this area.

\subsection{The Bilinear Maximal Functions}

The theory of the one dimensional Hilbert transform and maximal function 
are intimately related, hence it is natural to consider the bilinear 
maximal functions 
\begin{equation*}
\operatorname M _{\alpha }(f,g) := 
\sup _{t>0} (2t) ^{-1}\int _{-t} ^{t} \abs{ f(x-\alpha y)g(x-y)}\; dy\,.
\end{equation*}
For this operator, certain  bounds are immediately availible.  Namely, for 
conjugate indicies $ \tfrac1p_=\tfrac1 {p'}=1$, we have 
\begin{equation*}
\operatorname M _{\alpha } (f,g) \lesssim ( \operatorname M \abs{ f} ^{p}) ^{1/p}
(\operatorname M \abs{ g} ^{p'}) ^{1/p'}\,.
\end{equation*}
Thus, the known $ L^p$ inequalities for the maximal function are availible, showing 
that e.g.~the bilinear maximal function maps $ L^2\times L^2$ into weak $ L^1$. 

But, the bilinear Hilbert transform maps into certain spaces $ L^r$ for $ \tfrac23<r\le1$, 
so it is natural to ask if the bilinear maximal function does as well.  Indeed, this is true. 

\begin{theorem}\label{t.bmf} 
  For $1<p,q\le2$, if  $1\le 1/r=1/p+1/q<3/2$, and $\alpha\not\in\{0,1,\infty\}$, then 
 \begin{equation*}
\norm \operatorname M_\alpha (f,g).r.\lesssim{}\norm f.p.\norm g.q.
 \end{equation*}
\end{theorem}

The proof of the author \cite{MR1745019} begins by observing that the maximal function 
can be be bounded by truncations of a singular integral operator, with approriately 
chosen kernel.  Essentially, it is enough to bound the maximal truncations of the 
bilinear Hilbert transform given by 
\begin{equation*}
\sup _{k\in \mathbb Z} \ABs{\int  _{\abs{ y}\le2 ^{k}}f(x-\alpha y)g(x-y)\;\frac{dy}y }
\end{equation*}
These maximal operators obey the same inequalities in Theorem~\ref{t.bht}, and this is 
the main theorem in \cite {MR1745019}.  

The central point is to replace the Size Lemma above by a suitable maximal variant.  
This is can be done, but one must appeal to fundamental maximal inequality found by 
Bourgain \cite {MR1019960}.  The reader can also consult the paper of Thiele 
\cite {MR1846092} which presents the entire proof in the Walsh context, where many 
of the technical difficulties are minimized. 

Demeter, Tao and Thiele have revisited these issues. Strikingly, they found 
an argument which provides an $ \epsilon $ improvement over the trivial 
$ L^2\times L^2\mapsto \text {weak} L^1$ bound mentioned above.  Moreover, the 
argument uses arithmetic combinatorics.  See \cite{demeter-tao-thiele}.

\subsection{Multilinear Oscillatory Integrals}  \label{s.multi} 
Consider the bilinear oscillatory Integral 
 \begin{equation*}
\operatorname B_2 (f_1,f_2)(x):=\int f_1(x-y)f_2(x+y) \frac{ \operatorname e ^{2i y^2 } }y \; dy .
 \end{equation*}
This is seen to be a disguised form of a bilinear Hilbert transform.  Setting $g_j(x):=\operatorname  e ^{ix^2} f_j(x) $, one 
sees that 
 \begin{equation*}
\operatorname B_2(g_1,g_2)(x):=\operatorname e ^{2i x^2 } \int f_1(x-y)f_2(x+y) \; \frac{ dy }y  . 
 \end{equation*}
(Compare this to Ionescu's example mentioned at the end of \s.stein/.)   As it turns out, for a polynomial of 
any other degree, the integral above is bounded.  The proof demonstrates a multilinear variant of the van der Corput 
type inequality, of which \l.vdc/ is just one example.  

More generally,   M.~Christ, Li, Tao and Thiele \cite{cltt} define  multilinear functionals 
\begin{equation} \label{e.mulit-defn}
\Lambda_\lambda (f_1,f_2,\cdots,f_n)
= \int_{\mathbb R^m} e^{i\lambda P(x)} \prod_{j=1}^n f_j(\pi_j(x))\eta(x)\,dx
\end{equation}
where $\lambda\in \mathbb R $ is a parameter, $P:\mathbb R^m\mapsto \mathbb R$ is a  real-valued
polynomial, $m\ge 2$, and $\eta\in C^1_0(\mathbb R^m)$ is compactly supported. 
Each $\pi_j$ denotes the orthogonal projection from $\mathbb R^m$
to a linear subspace $V_j\subset\mathbb R^m$ of any dimension $\kappa\le m-1$,
and $f_j:V_j\to \mathbb C$ is always assumed to be locally integrable
with respect to Lebesgue measure on $V_j$.  

Notice that by taking $n=3 $, and taking projections $\pi_j \mid \mathbb R^2\mapsto V_j $  where 
 \begin{equation} \label{e.V-ex}
\begin{split} 
V_1=\{ (x,x)\mid x\in \mathbb R\},\qquad V_2=\{(x,-x)\mid x\in \mathbb R\},
\\ V_3=\{(x,0)\mid x\in \mathbb R\}, 
\end{split}
 \end{equation}
we can recover for instance a bilinear Hilbert transform.

And they say that a polynomial $P $ has \emph{ a power decay property} if there is a $\delta>0 $, so that 
for all $f_j \in L^ \infty(V_j ) $, we have the estimate 
 \begin{equation*}
\abs{ \Lambda(f_1,\ldots, f_n)} {}\lesssim{} (1+\abs{\lambda}) ^{-\delta }\prod _{j=1}^ n \norm f_j .\infty. 
 \end{equation*}
From this estimate, a range of power decay estimates hold in all relevant products of $L^p $ spaces. 
This should be compared to \l.vdc/ and in particular (\ref{e.vdc})  below. 

Clearly, there are obstructions to a power decay property, and this obstruction can be formalized in a definition.  
A polynomial $P$ is said to be\emph{ degenerate (relative to $\{V_j\}$)} 
if there exist polynomials $p_j:V_j\to\mathbb R$
such that $P=\sum_{j=1}^n p_j\circ\pi_j$.
Otherwise $P$ is nondegenerate.
In the case $n=0$, where the collection of subspaces $\{V_j\}$ is empty,
$P$ is considered to be nondegenerate if and only if it is nonconstant.
And in the example (\ref{e.V-ex}) , we see that $P(y)=2x^2+2y^2=(x+y)^2+(x-y)^2 $ is degenerate.  

It is natural to conjecture that non degeneracy is  sufficient for  a power decay property.  This is  verified 
in a wide range of special cases in the paper Christ, Li, Tao, and Thiele \cite{cltt}, by a range of interesting 
techniques.  It is of interest to determine if the natural conjecture here is indeed correct.

\subsection{ Hilbert Transform on Smooth Families of Lines}

This quetion has its beginnings in the Besicovitch  set, which we already mentioned in connection to spherical summation of Fourier series.  
One may construct Besicovitch sets with these properties.  For choices of $0<\epsilon,\alpha<1$, 
there is a Besicovitch set $K$ in the square $[0,4]^2$ say, for which $K$ has measure at most
$\epsilon$, and there is a function $g\mid \mathbb R^2\to \mathbb T$, so that  for a set of $x$'s in $[0,4]^2$ of measure $\gtrsim1$,
$K\cap \{x+tv(x)\mid t\in\mathbb R\}$ contains a line segment of length one, and $v$ is H\"older continuous of order $\alpha$.  

One can ask if the H\"older continuity condition is sharp.  A beautiful formulation of a conjecture in this direction is attributed  to A.~Zygmund.

\begin{conjecture}  Let $v\mid \mathbb R^2\to\mathbb T$ be H\"older   
continuous (of order $1$).  Then for all square integrable functions $f$, 
$$
f(x)=\lim_{t\to0}(2t)^{-1}\int_{-t}^t f(x-uv(x))\; du\qquad \text{a.e}(x)
$$
\end{conjecture}

This is a differentiablity question, on a choice of lines specified by $v$.  The only stipulation is that $v$ is H\"older continuous.
This is only known under more stringent conditions on $v$, such as analytic due to E.M.~Stein  \cite{steinICM}, or real
 analytic due to Bourgain \cite{bourgain}.   There is a partial result due to N.~Katz  \cite{MR2004e:42030} (also see \cite{nets1}) 
 that demonstrates 
at worst ``log log''
blowup assuming the H\"older continuity of $v$.  
The question is  open, even if one assumes that $v\in C^{1000}$.  

The difficulty in this problem arises from those points at which the gradient of $v$ is degenerate; 
assumptions such as analyticity certainly control such degeneracies.  

E.M.~Stein \cite{steinICM} posed the Hilbert transform variant, namely defining 
$$
H_v f(x):=\text{p.v.}\int_{-1}^1 f(x-yv(x))\frac{dy}y, 
$$
is it the case that there is a constant $c$ so that if $\norm v.\text{H\"ol}.<c$, then   
$H_v$ maps $L^2(\mathbb R^2)$ into itself.   A curious fact about this question is that this inequality,
if known, implies Carleson's theorem for one dimensional Fourier series.  

To see this, observe that the symbol for the transform is $\psi(\xi\cdot v(x))$, where $\psi$ is the Fourier 
transform of $y^{-1}\ind {\{\abs y<1\}}  $.   Suppose the vector field is of the form 
$v(x)=(1,\nu(x_1))$ where we need only assume that $\nu$ is Holder continuous of norm $1$ say, and consider
the trace of the symbol on the line $\xi_2=-N$.  Then, the symbol is $\psi((\xi_1,N)\cdot(1,\nu(x_1))=\psi(\xi_1-N\nu(x_1))$.  We conclude that thi
s symbol 
defines a bounded linear operator on $L^2(\mathbb R)$, with bound that is independent of $N$.  That is, for {\it any} 
Lipschitz function $\nu(x_1)$,  and {\it any $N>1$} the symbol $\psi(\xi_1-N\nu(x_1))$ is the symbol of a bounded 
linear operator on $L^2(\mathbb R)$
.  By varying $N$ and $\nu$, we may replace $N\nu(x_1)$ by an arbitrary measurable function. This is the substance of
Carleson's theorem. 

But the implication is entirely one way:  A positive answer to the family of lines question seems to require techniques quite
a bit more sophisticated than those that imply Carleson's theorem.   Recently Lacey and Li \cites{laceyli,
math.CA/0310345} have been able to 
obtain a partial answer, 
assuming only that the vector field has $1+\epsilon$ derivatives. 

 \begin{theorem}\label{t.laceyli}  Assume that $v\in C^{1+\epsilon}$ for some $\epsilon>0$.  Then the operator $H_v$ is bounded on  $L^2(\mathbb R^2)$.  
The norm of the operator is at most
 \begin{equation*} 
\norm H_v.2.\lesssim{} [1+ {\log_+ \norm v.C^{1+\epsilon}.}]^{2}.
 \end{equation*}
 \end{theorem}

\subsection{Schr\"odinger Operators, Scattering Transform}

There is a beautiful line of investigation relating Schr\"odinger equations in one dimension to 
aspects of the Fourier transform, and in particular, 
Carleson's theorem. There is a further connection to scattering transforms and 
 nonlinear Fourier analysis.  All in all, these topics are extremely broad, with several different 
 sets of motivations, and a long list of contributors.  
 
 We concentrate on a succinct way to see the connection to Carleson's theorem, an observation made 
 explicitly  by M.~Christ and A.~Kiselev \cites{MR2001k:34099,MR1952927},  also see C.~Remling \cite{rem}.  
 The basic object is a time independent 
 Schr\"odinger operator on the real line,   
  \begin{equation*} 
 H=-\frac {d^2}{dx^2}+V
  \end{equation*}
 where $V$ is an appropriate potential on the real line.  The idea is that if $V$ is small, in some 
 specific senses, then the spectrum of $H$ should resemble that of $-\frac {d^2}{dx^2}$.  In particular 
 eigenfunctions should be perturbations of the exponentials. 
 
 Standard examples show that one should seek to show that for almost all $\lambda$, the eigenfunctions of energy $\lambda$,
 that is the solutions to 
  \begin{equation*}
 (H-\lambda^2 I)
  \end{equation*}
 are bounded perturbations of $e^{\pm i\lambda x}$.  
 
 Seeking such an eigenfunction, one can formally write 
  \begin{equation*}
 u(x)=e^{i\lambda x}+\frac 1{i\lambda}\int_x^\infty \sin(\lambda (x-y)) V(y)u(y)\; dy
  \end{equation*}
 Iterating this formula, again formally, one has 
  \begin{gather}
 u(x)={}e^{i\lambda x}  \nonumber\\
 {}+\frac 1{i\lambda}\int_x^\infty \sin(\lambda (x-y)) V(y) e^{i\lambda y} \; dy
 \label{e.single} \\  \label{e.double}
{}+\frac 1{(i\lambda)^2}\int\int_{x\le{}y_1\le{}y_2}\sin (x-y_1)\sin(y_1-y_2)V(y_1)V(y_2) u(y_2)\; dy_1\,dy_2 .
 \end{gather}

Observe that (\ref{e.single})  no longer contains $u$, and  is a linear combination of 
 \begin{gather}
\label{e.singcar} 
e^{i\lambda x}\int_x^\infty e^{2i\lambda y} V(y)\; dy 
\\ 
\label{e.singl1} 
e^{i\lambda x}\int_x^\infty V(y)\; dy 
 \end{gather}
One seeks estimates of these in the mixed norm space of say, $L^2_\lambda L^\infty_x$.  From such estimates, one 
deduces that for almost all $\lambda$, there is an eigenfunction with is a perturbation of $e^{i\lambda x}$. 

Concerning (\ref{e.singcar}) , notice that if $V\in L^2$, we can, by Plancherel, regard $V$ as $\widehat f$, for some 
$f\in L^2$. The desired estimate is then a consequence of Carleson's theorem. This is indicative of the distinguished 
role that $L^2$ plays in this subject.  Also of the intertwining of the roles of frequency and time that occur in the subject. 

Concerning (\ref{e.singl1}) , unless $V\in L^1$, there is no reasonable interpretation that can be placed on this term.  In practice, 
a different approach than the one given here  must be adopted.  

If one continues the expansion in (\ref{e.double}) , one gets a bilinear operator with features that resemble both the Carleson 
operator, and the bilinear Hilbert transform.  See the papers by C.~Muscalu,  T.~Tao, and C.~Thiele \cites{MR1952931,biest1,biest2,counter}.

We refer the reader to these papers by M.~Christ and A.~Kiselev \cites{MR2001k:34099,MR1952927,MR2002e:34143,MR2001i:47054}.  
For a survey of this subject, see  M.~Christ and A.~Kiselev \cite{survey}.
 The reader should also consult the ongoing investigations of C.~Muscalu, T.~Tao, and C.~Thiele \cite{MR1950785}. 
This  paper begins with an interesting summary of the perspective of the nonlinear Fourier transform.

\bigskip
 {
Michael T. Lacey \hfill\break
School of Mathematics\hfill\break
Georgia Institute of Technology\hfill\break
Atlanta GA 30332\hfill\break
\smallskip
\tt lacey@math.gatech.edu\hfill\break
\tt http://www.math.gatech.edu/\~{}lacey }
 
\end{document}